\documentclass[12pt,oneside,reqno]{amsart}
\usepackage{mathrsfs}
\usepackage{graphics}
\usepackage{amssymb}
\usepackage{enumerate}
\newcommand{\dif}{\mathrm{d}}

\newcommand{\be}{\begin{eqnarray}}
\newcommand{\ee}{\end{eqnarray}}
\newcommand{\ce}{\begin{eqnarray*}}
\newcommand{\de}{\end{eqnarray*}}
\newtheorem{theorem}{Theorem}[section]
\newtheorem{lemma}[theorem]{Lemma}
\newtheorem{remark}[theorem]{Remark}
\newtheorem{definition}[theorem]{Definition}
\newtheorem{proposition}[theorem]{Proposition}
\newtheorem{Example}[theorem]{Example}
\newtheorem{corollary}[theorem]{Corollary}
\def\e{\varepsilon}
\def\s{\sigma}

\def\a{\alpha}

\def\b{\beta}
\def\d{\delta}
\def\p{\partial}
\def\g{\gamma}

\def\[{{\Big[}}
\def\]{{\Big]}}
\def\<{{\langle}}
\def\>{{\rangle}}
\def\({{\Big(}}
\def\){{\Big)}}

\def\no{\nonumber}
\def\bt{\begin{theorem}}
\def\et{\end{theorem}}
\def\bl{\begin{lemma}}
\def\el{\end{lemma}}
\def\br{\begin{remark}}
\def\er{\end{remark}}
\def\bx{\begin{Example}}
\def\ex{\end{Example}}
\def\bd{\begin{definition}}
\def\ed{\end{definition}}
\def\bp{\begin{proposition}}
\def\ep{\end{proposition}}
\def\bc{\begin{corollary}}
\def\ec{\end{corollary}}

\def\cB{{\mathcal B}}

\def\cH{{\mathcal H}}

\def\cL{{\mathcal L}}
\def\cM{{\mathcal M}}
\def\cN{{\mathcal N}}

\def\cP{{\mathcal P}}

\def\mC{{\mathbb C}}

\def\mE{{\mathbb E}}

\def\mP{{\mathbb P}}

\def\mR{{\mathbb R}}

\def\mW{{\mathbb W}}

\def\sB{{\mathscr B}}

\def\sF{{\mathscr F}}

\def\sL{{\mathscr L}}

\def\sP{{\mathscr P}}

\def\geq{\geqslant}
\def\leq{\leqslant}

\begin{document}

\allowdisplaybreaks

\title{Weak approximation of nonlinear filtering for multiscale McKean-Vlasov stochastic systems*}

\author{Huijie Qiao and Wanlin Wei }

\dedicatory{School of Mathematics,
Southeast University,\\
Nanjing, Jiangsu 211189, P.R.China}



\thanks{*This work was partly supported by NSF of China (No.12071071).}

\thanks{Corresponding author: Huijie Qiao, hjqiaogean@seu.edu.cn}


\date{}

\begin{abstract}
The work concerns the nonlinear filtering problem for a class of multiscale McKean-Vlasov stochastic systems. First of all, by a Poisson equation we prove that the solution of the slow part for a multiscale system weakly converges to the solution of the average equation. Then we define nonlinear filtering of the origin multiscale system and the average equation, and again through the same Poisson equation show the weak approximation between nonlinear filtering of the slow part for the origin multiscale system and that of the average equation.
\end{abstract}

\maketitle \rm

\section{Introduction}

Let $(\Omega,\sF,\{\sF_{t}\}_{t\in[0,T]},\mP)$ be a complete filtered probability space. $B, W$ are $n$-dimensional and $m$-dimensional standard Brownian motions defined on $(\Omega,\sF,\{\sF_{t}\}_{t\in[0,T]},\mP)$, respectively. Moreover, $B$ and  $W$ are mutually independent. Consider the following system on $\mR^{n} \times \mR^{m}$:
\be\left\{\begin{array}{l}
\dif X_{t}^{\e}=b_{1}(X_{t}^{\e}, \sL^{\mP}_{X_{t}^{\e}}, Z_{t}^{\e})\dif t+\s_{1}(X_{t}^{\e}, \sL^{\mP}_{X_{t}^{\e}}, Z_{t}^{\e})\dif B_{t},\\
X_{0}^{\e}=x_0,\quad  0\leq t\leq T,\\
\dif Z_{t}^{\e}=\frac{1}{\e}b_{2}(\sL^{\mP}_{X_{t}^{\e}},Z_{t}^{\e})\dif t+\frac{1}{\sqrt{\e}}\s_{2}(\sL^{\mP}_{X_{t}^{\e}},Z_{t}^{\e})\dif W_{t},\\
Z_{0}^{\e}=z_0,\quad  0\leq t\leq T,
\end{array}
\right.
\label{slfa}
\ee
where these mappings $b_{1}:\mR^{n}\times\cP_2(\mR^{n})\times\mR^{m}\rightarrow\mR^{n}$, $\s_{1} :\mR^{n}\times\cP_2(\mR^{n})\times\mR^{m}\rightarrow\mR^{n\times n}$, $b_{2} :\cP_2(\mR^{n})\times\mR^{m}\rightarrow\mR^{m}$, $\s_{2} :\cP_2(\mR^{n})\times\mR^{m}\rightarrow\mR^{m\times m}$ are all Borel measurable, $\sL^{\mP}_{X_{t}^{\e}}$ denotes the distribution of $X_{t}^{\e}$ under the probability measure $\mP$, and $\cP_{2}(\mR^n)$ denotes  the collection of probability measures on $\mR^n$ with finite second moments. 

If $b_1, \s_1, b_2, \s_2$ don't depend on the distribution $\sL^{\mP}_{X_{t}^{\e}}$, and $b_2, \s_2$ also depend on the state $X^\e_t$, the system (\ref{slfa}) is usually called a multiscale or slow-fast system. And $X^\e_{\cdot}, Z_{\cdot}^{\e}$ are called the slow and fast parts, respectively. These systems are often used to model phenomena from many fields, such as biology, chemistry and engineering. Moreover, there are many results about them. We only mention the average principle for them. In \cite{rzk}, Khasminskii first studied the average principle. Later, many authors generalized this result (See e.g. \cite{LRSX, LD1, lh, pv1, pv2, pv3} and the references therein). In particular, Pardoux and Veretennikov \cite{pv1,pv2,pv3} systematically investigated the approximation between $X^\e_t$ and the solution of the average equation through Poisson equations. 

If $b_1, \s_1, b_2, \s_2$ depend on distributions, the system (\ref{slfa}) is generally called a multiscale McKean-Vlasov stochastic system. Since McKean-Vlasov stochastic differential equations (SDEs for short) are more complex than SDEs, there are only a few results about the average principle. When $\s_{1}$ doesn't depend on $Z_{t}^{\e}$ and $b_2, \s_2$ also depend on the state $X^\e_t$, R\"{o}ckner, Sun and Xie \cite{SunX} showed that the slow part $X^\e_{\cdot}$ of the system (\ref{slfa}) converges to an average system in the $L^{2}$ sense. Recently, Xu et al. \cite{XLLM} extended the result in \cite{SunX} to the case where $b_2, \s_2$ depend on $(X^\e_t,\sL^{\mP}_{X_{t}^{\e}},Z^\e_t,\sL^{\mP}_{Z_{t}^{\e}})$. Very recently, we \cite{Qiao2} considered a more general system than that in \cite{XLLM} and established the $L^{p}$ ($p\geq 16$) convergence. In this paper, we observe the system (\ref{slfa}), that is, $\s_{1}$ does depend on $Z_{t}^{\e}$. The price to pay is that $b_2, \s_2$ don't depend on the state $X^\e_t$ and the distribution $\sL^{\mP}_{Z_{t}^{\e}}$.

 Next, we take an observation process $Y_{t}^{\e}$, i.e.
\be
Y_{t}^{\e}=V_{t}+\int_{0}^{t}h(X_{s}^{\e},\sL^{\mP}_{X_{s}^{\e}},Z_{s}^{\e})\dif s,
\label{byte}
\ee
where $V$ is a $l$-dimensional Brownian motion independent of $B, W$, and $h: \mR^n\times\cP_2(\mR^n)\times\mR^m\rightarrow \mR^l$ is Borel measurable. Thus, the system (\ref{slfa}) and the process (\ref{byte}) form a signal-observation system. And the nonlinear filtering problem for this signal-observation system is to estimate and predict $(X^\e_{\cdot},\sL^{\mP}_{X_{\cdot}^{\e}},Z^\e_{\cdot})$ based on the information of $Y^\e_{\cdot}$. This problem arises in many fields, such as stochastic control, financial modeling, speech and image processing, and Bayesian networks. Furthermore, if $b_1, \s_1, b_2, \s_2, h$ don't depend on the distribution $\sL^{\mP}_{X_{t}^{\e}}$, this problem has been widely studied (See \cite{bnp, Goggin1, Goggin2, Imkeller, lh, Qiao1, Qiao3} and references therein). If $\s_{1}, h$ don't depend on $Z_{t}^{\e}$ and $b_2, \s_2$ also depend on the state $X^\e_t$, we \cite{Qiao2} showed the convergence of the nonlinear filtering for $(X^\e_t,\sL^{\mP}_{X_{t}^{\e}})$ in the $L^{p}$ ($p\geq 1$) sense. In this paper, we require that $\s_{1}, h$ depend on $Z_{t}^{\e}$ and obtain the weak convergence.

As a whole, our contribution are two-folded:

$\bullet$ By a Poisson equation we prove that the slow part of the original system weakly converges to the average system.

$\bullet$ By the same Poisson equation we establish that the nonlinear filtering of the slow part weakly converges to that of the average system.

It is worthwhile to mentioning our results. In \cite{lwx}, three authors studied the following multiscale McKean-Vlasov system:
\be\left\{\begin{array}{l}
\dif \hat{X}_{t}^{\e}=\hat{b}_{1}(\hat{X}_{t}^{\e}, \sL^{\mP}_{\hat{X}_{t}^{\e}}, \hat{Z}_{t}^{\e},\sL^{\mP}_{\hat{Z}_{t}^{\e}})\dif t+\hat{\s}_{1}(\hat{X}_{t}^{\e}, \sL^{\mP}_{\hat{X}_{t}^{\e}}, \hat{Z}_{t}^{\e},\sL^{\mP}_{\hat{Z}_{t}^{\e}})\dif B_{t},\\
\hat{X}_{0}^{\e}=x_0,\quad  0\leq t\leq T,\\
\dif \hat{Z}_{t}^{\e}=\frac{1}{\e}\hat{b}_{2}(\sL^{\mP}_{\hat{X}_{t}^{\e}},\hat{Z}_{t}^{\e},\sL^{\mP}_{\hat{Z}_{t}^{\e}})\dif t+\frac{1}{\sqrt{\e}}\hat{\s}_{2}(\sL^{\mP}_{\hat{X}_{t}^{\e}},\hat{Z}_{t}^{\e},\sL^{\mP}_{\hat{Z}_{t}^{\e}})\dif W_{t},\\
\hat{Z}_{0}^{\e}=z_0,\quad  0\leq t\leq T,
\end{array}
\right.
\label{slfa1}
\ee
where these mappings $\hat{b}_{1}:\mR^{n}\times\cP_2(\mR^{n})\times\mR^{m}\times\cP_2(\mR^{m})\rightarrow\mR^{n}$, $\hat{\s}_{1} :\mR^{n}\times\cP_2(\mR^{n})\times\mR^{m}\times\cP_2(\mR^{n})\rightarrow\mR^{n\times n}$, $\hat{b}_{2} :\cP_2(\mR^{n})\times\mR^{m}\times\cP_2(\mR^{m})\rightarrow\mR^{m}$, $\hat{\s}_{2} :\cP_2(\mR^{n})\times\mR^{m}\times\cP_2(\mR^{m})\rightarrow\mR^{m\times m}$ are all Borel measurable. There they defined 
\be
\bar{\hat{b}}_1(x,\mu):=\int_{\mR^m}\hat{b}_1(x,\mu,z,\hat{\nu}^\mu)\hat{\nu}^\mu(\dif z), \quad \bar{\hat{\s}}_1\bar{\hat{\s}}^*_1(x,\mu):=\int_{\mR^m}(\hat{\s}_1\hat{\s}^*_1)(x,\mu,z,\hat{\nu}^\mu)\hat{\nu}^\mu(\dif z),
\label{bbar2}
\ee
where $\hat{\nu}^\mu$ is the unique invariant probability measure for the frozen equation (See Section \ref{mainresu}), constructed the average equation (\ref{Eq3})  and presented a similar average principle under strong conditions. It is obvious that the system (\ref{slfa1}) is more general than the system (\ref{slfa}). However, comparing (\ref{bbar2}) with (\ref{bbar}), we think that our definitions of $\bar{b}_1, \bar{\s}_1$ are more natural. Moreover, our proof of the average principle is more succinct. Therefore, we don't simply seek generality. Besides, Beeson et al. \cite{bnp} considered the nonlinear filtering problem of the following system:
\be\left\{\begin{array}{l}
\dif \check{X}_{t}^{\e}=\check{b}_{1}(\check{X}_{t}^{\e}, \check{Z}_{t}^{\e})\dif t+\check{\s}_{1}(\check{X}_{t}^{\e}, \check{Z}_{t}^{\e})\dif B_{t},\\
\check{X}_{0}^{\e}=\check{X}_0,\quad  0\leq t\leq T,\\
\dif \check{Z}_{t}^{\e}=\frac{1}{\e}\check{b}_{2}(\check{X}_{t}^{\e},\check{Z}_{t}^{\e})\dif t+\frac{1}{\sqrt{\e}}\check{\s}_{2}(\check{X}_{t}^{\e},\check{Z}_{t}^{\e})\dif W_{t},\\
\check{Z}_{0}^{\e}=\check{Z}_0,\quad  0\leq t\leq T,\\
\check{Y}_{t}^{\e}=V_{t}+\int_{0}^{t}\check{h}(\check{X}_{s}^{\e},\check{Z}_{s}^{\e})\dif s,
\end{array}
\right.
\label{slfa2}
\ee
where these mappings $\check{b}_{1}:\mR^{n}\times\mR^{m}\rightarrow\mR^{n}$, $\check{\s}_{1} :\mR^{n}\times\mR^{m}\rightarrow\mR^{n\times n}$, $\check{b}_{2} :\mR^{n}\times\mR^{m}\rightarrow\mR^{m}$, $\check{\s}_{2} :\mR^{n}\times\mR^{m}\rightarrow\mR^{m\times m}$, and $\check{h}: \mR^n\times\mR^m\rightarrow \mR^l$ are all Borel measurable. They proved that the nonlinear filtering of the slow part $X_t^\e$ converges in probability to ``the nonlinear filtering" of the solution for the average equation. Here, for the signal-observation system (\ref{slfa})+(\ref{byte}) all the coefficients depend on the distribution $\sL^{\mP}_{X_{t}^{\e}}$. Therefore, our model is more general in some sense.

Finally, let us describe our motivation of this paper. McKean-Vlasov SDEs appear in the analysis of interacting particles in mathematical physics where so-called mean field type interaction is shown to be modeled by such a nonlinear process. Thus, in order to simulate an interacting system with a fast part and a slow part, we need to use a multiscale McKean-Vlasov stochastic system. Since the data of the original system is large, the data of the multiscale system is even larger. Thus, it becomes a problem to estimate and predict such systems. Hence, nonlinear filtering theory for multiscale McKean-Vlasov systems has direct applications in the estimation problems of mathematical physics and nonlinear mean field game theory and it is this which motivates one to study nonlinear filtering theory for systems with multiscale McKean-Vlasov dynamics.

The paper proceeds as follows. In Section \ref{pre}, we introduce notations and the definition of $L$-derivatives. Then we state main results in Section \ref{mainresu}. The proofs of two main theorems are placed in Section \ref{prooaver} and \ref{proofilt}, respectively. We present an example to explain our results in Section \ref{exam}. Finally, in Section \ref{app}, we complete the proof of an inequality.

The following convention will be used throughout the paper: $C$ with or without indices will denote different positive constants whose values may change from one place to another.

\section{Preliminary}\label{pre}

In this section, we will recall some notations, and the definition of $L$-derivatives for functions on $\cP_2(\mR^n)$ and list all the assumptions.

\subsection{Notations}\label{nn}

In this subsection, we introduce some notations used in the sequel.

Let $|\cdot|, \|\cdot\|$ be the norm of a vector and a matrix, respectively. Let $\langle\cdot,\cdot\rangle$ be the inner product of vectors on
$\mR^n$. Let $A^{*}$ be the transpose of the matrix $A$.

Let $\cB_b(\mR^{n})$ be the set of all bounded Borel measurable functions on $\mR^n$. Let $C(\mR^n)$ be the set of all  functions which are continuous on $\mR^n$. $C^{k}(\mR^{n})$ represents the collection of all functions which are continuous differentiable up to $k$-order. 

Let $\sB(\mR^n)$ be the Borel $\sigma$-field on $\mR^n$. Let $\cP(\mR^n)$ be the collection of all probability measures on $\sB(\mR^n)$ with the usual topology of weak convergence. Let $\cP_{2}(\mR^n)$ denote  the collection of probability measures on $\sB(\mR^n)$ satisfying:
$$
\|\mu\|^{2}:=\int_{\mR^n}|x|^{2}\mu(dx)<\infty.
$$
It is known that $\cP_2(\mR^n)$ is a Polish space endowed with the $L^2$-Wasserstein distance defined by
$$
\mathbb{W}_2(\mu,\nu):= \inf\limits_{\pi\in\Psi(\mu,\nu)}\left(\int_{\mathbb{R}^n\times\mathbb{R}^n}|x-y|^{2}\pi(\dif x,\dif y)\right)^{\frac{1}{2}}, \quad \mu , \nu\in \cP_2(\mR^n),
$$
where $\Psi(\mu,\nu)$ is the set of all couplings $\pi$ with marginal distributions $\mu$ and $\nu$. Moreover, if $\xi,\zeta$ are two random variables with distributions $\sL_\xi, \sL_\zeta$ under $\mP$, respectively,
$$
\mathbb{W}_2(\sL_\xi, \sL_\zeta)\leq (\mE|\xi-\zeta|^2)^{\frac{1}{2}},
$$
where $\mE$ stands for the expectation with respect to $\mP$.

\subsection{$L$-derivatives for functions on $\cP_2(\mR^n)$}\label{lde} In this subsection we recall the definition of $L$-derivatives for functions on $\cP_2(\mR^n)$ (c.f. \cite{rw}). 

Let $I$ be the identity map on $\mR^n$. For $\mu\in\cP_2(\mR^n)$ and $\phi\in L^2(\mR^n, \sB(\mR^n), \mu;\mR^n)$, $<\mu,\phi>:=\int_{\mR^n}\phi(x)\mu(\dif x)$. Moreover, by simple calculation, it holds that $\mu\circ(I+\phi)^{-1}\in\cP_2(\mR^n)$.

\bd\label{lderi}
(i) A function $f: \cP_2(\mR^n)\mapsto\mR$ is called L-differentiable at $\mu\in\cP_2(\mR^n)$, if the functional
$$
L^2(\mR^n, \sB(\mR^n), \mu;\mR^n)\ni\phi\mapsto f(\mu\circ(I+\phi)^{-1})
$$
is Fr\'echet differentiable at $\phi=0$; that is, there exists a unique $\gamma\in L^2(\mR^n, \sB(\mR^n), \mu;\mR^n)$ such that
$$
\lim\limits_{<\mu,|\phi|^2>\rightarrow 0}\frac{f(\mu\circ(I+\phi)^{-1})-f(\mu)-<\mu,\gamma\cdot\phi>}{\sqrt{<\mu,|\phi|^2>}}=0.
$$
In the case, we denote $\partial_{\mu}f(\mu)=\gamma$ and call it the $L$-derivative of $f$ at $\mu$.

(ii) A function $f: \cP_2(\mR^n)\mapsto\mR$ is called $L$-differentiable on $\cP_2(\mR^n)$ if $L$-derivative $\partial_{\mu}f(\mu)$ exists for all $\mu\in\cP_2(\mR^n)$.
\ed

\bd\label{space1}
 The function $f$ is said to be in $C^{(k_1,k_2)}(\cP_2(\mR^n))$, if $f$ is continuous $L$-differentiable to $k_1$-order, and for any $\mu\in\cP_2(\mR^n)$, $(y_1,y_2,\cdots,y_{k_1})\mapsto\partial^{k_1}_\mu f(\mu)(y_1,y_2,\cdots,y_{k_1})$ is in $C^{k_2}(\mR^{n k_1})$. 
\ed

\bd\label{space2}
The function $F: \mR^n\times\cP_2(\mR^n)\mapsto\mR$ is said to be in $C^{k_0,(k_1,k_2)}(\mR^n\times\cP_2(\mR^n),\mR)$, if for any $\mu\in\cP_2(\mR^n)$, $x\mapsto F(x,\mu)$ is in $C^{k_0}(\mR^n)$ and for any $x\in\mR^n$, $\mu\mapsto F(x,\mu)$ is in $C^{(k_1,k_2)}(\cP_2(\mR^n))$. If $F\in C^{k_0,(k_1,k_2)}(\mR^n\times\cP_2(\mR^n),\mR)$, and itself and all the derivatives are uniformly bounded and jointly continuous in the corresponding variable family, we say that $F\in C_b^{k_0,(k_1,k_2)}(\mR^n\times\cP_2(\mR^n),\mR)$.
\ed

\bd\label{space3}
The function $F: \mR^n\times\cP_2(\mR^n)\mapsto\mR$ is said to be in $\mC_b^{4,(2,2)}(\mR^n\times\cP_2(\mR^n),\mR)$, if $F\in C_b^{4,(2,2)}(\mR^n\times\cP_2(\mR^n),\mR)$ and $\p_x\p_\mu F(x,\mu)(y), \p_x\p_y\p_\mu F(x,\mu)(y), \p_\mu\p_x F(x,\mu)(y), \\\p_\mu\p_{xx} F(x,\mu)(y), \p_{xx}\p_\mu F(x,\mu)(y), \p_{xx}\p_y\p_\mu F(x,\mu)(y), \p_y\p_\mu\p_x F(x,\mu)(y), \p_y\p_\mu\p_{xx} F(x,\mu)(y) $ exist and are uniformly bounded and jointly continuous in the corresponding variable family.
\ed

\bd\label{space4}
The function $\Psi: \mR^n\times\cP_2(\mR^n)\times\mR^m\mapsto\mR$ is said to be in $C^{k_0,(k_1,k_2),k_0}(\mR^n\times\cP_2(\mR^n)\times\mR^m)$, if for $(x,z)\in\mR^n\times\mR^{m}$,$\Psi(x,\cdot,z)\in C^{(k_1,k_2)}(\cP_2(\mR^n))$ and for $\mu\in\cP_2(\mR^n)$,$\Psi(\cdot,\mu,\cdot)\in C^{k_0}(\mR^n\times\mR^{m})$. If $\Psi\in C^{k_0,(k_1,k_2),k_0}(\mR^n\times\cP_2(\mR^n)\times\mR^m)$, and itself and all its derivatives are bounded and jointly continuous in the corresponding variable family, we say $\Psi\in C_b^{k_0,(k_1,k_2),k_0}(\mR^n\times\cP_2(\mR^n)\times\mR^m)$.
\ed

\subsection{Assumptions}

In the subsection, we list all the assumptions used in the sequel.

\begin{enumerate}[$(\mathbf{H}^1_{b_{1}, \s_{1}})$]
\item
There exists a constant $L_{b_{1}, \s_{1}}>0$ such that for $x_{i}\in\mR^n$, $\mu_{i}\in\cP_{2}(\mR^n)$, $z_{i}\in\mR^m$, $i=1, 2$,
\ce
&&|b_{1}(x_{1},\mu_{1},z_{1})-b_{1}(x_{2},\mu_{2},z_{2})|^{2}+\|\s_{1}(x_{1},\mu_{1},z_{1})-\s_{1}(x_{2},\mu_{2},z_{2})\|^{2}\\
&\leq& L_{b_{1},\s_{1}}\(|x_{1}-x_{2}|^{2}+\mW_{2}(\mu_{1},\mu_{2})^2+|z_{1}-z_{2}|^{2}\).
\de
\end{enumerate}
\begin{enumerate}[$(\mathbf{H}^2_{\s_{1}})$]
\item There exists a constant $l>0$ such that for $x\in\mR^n$, $\mu\in\cP_{2}(\mR^n)$, $z\in\mR^m$, $\eta\in\mR^n$,
$$
\<\s_{1}(x,\mu, z)\eta, \eta\>\geq l|\eta|^{2}.
$$
\end{enumerate}
\begin{enumerate}[$(\mathbf{H}^3_{b_{1}, \s_{1}})$]
\item
$\partial_x b_1(x,\mu,z), \partial_{\mu}b_1(x,\mu,z)(y), \partial_z b_1(x,\mu,z), \partial_{xx} b_1(x,\mu,z), \partial_{zz} b_1(x,\mu,z), \partial_y\partial_{\mu}b_1(x,\mu,z)(y)$ exist for any $(x,\mu,y,z)\in\mR^n\times\cP_{2}(\mR^n)\times\mR^n\times\mR^m$ and $\partial_{xx} b_1(x,\mu,z),\partial_{zz} b_1(x,\mu,z)$ are uniformly bounded. And there exists a constant $\g_1\in(0,1]$ such that for $z_{i}\in\mR^m$, $i=1, 2$,
\ce
&&\sup\limits_{x\in\mR^n, \mu\in\cP_{2}(\mR^n)}\|\partial_{x}b_1(x,\mu,z_1)-\partial_{x}b_1(x,\mu,z_2)\|\leq C|z_1-z_2|^{\g_1},\\
&&\sup\limits_{x\in\mR^n, \mu\in\cP_{2}(\mR^n)}\|\partial_{\mu}b_1(x,\mu,z_1)-\partial_{\mu}b_1(x,\mu,z_2)\|_{L^2(\mu)}\leq C|z_1-z_2|^{\g_1},\\
&&\sup\limits_{x\in\mR^n, \mu\in\cP_{2}(\mR^n)}\|\partial_{z}b_1(x,\mu,z_1)-\partial_{z}b_1(x,\mu,z_2)\|\leq C|z_1-z_2|^{\g_1},\\
&&\sup\limits_{x\in\mR^n, \mu\in\cP_{2}(\mR^n)}\|\partial_{xx} b_1(x,\mu,z_1)-\partial_{xx} b_1(x,\mu,z_2)\|\leq C|z_1-z_2|^{\g_1},\\
&&\sup\limits_{x\in\mR^n, \mu\in\cP_{2}(\mR^n)}\|\partial_{zz} b_1(x,\mu,z_1)-\partial_{zz} b_1(x,\mu,z_2)\|\leq C|z_1-z_2|^{\g_1},\\
&&\sup\limits_{x\in\mR^n, \mu\in\cP_{2}(\mR^n)}\|\partial_y\partial_{\mu}b_1(x,\mu,z_1)-\partial_y\partial_{\mu}b_1(x,\mu,z_2)\|_{L^2(\mu)}\leq C|z_1-z_2|^{\g_1}.
\de
Moreover, the above conditions also hold for $\s_1\s_1^*$. 
\end{enumerate}
\begin{enumerate}[$(\mathbf{H}^1_{b_{2}, \s_{2}})$]
\item
There exists a constant $L_{b_{2}, \s_{2}}>0$ such that for $\mu_{i}\in\cP_{2}(\mR^n)$, $z_{i}\in\mR^m$, $i=1, 2$,
\ce
|b_{2}(\mu_{1},z_{1})-b_{2}(\mu_{2},z_{2})|^{2}+\|\s_{2}(\mu_{1},z_{1})-\s_{2}(\mu_{2},z_{2})\|^{2}\leq L_{b_{2}, \s_{2}}\(\mW_{2}(\mu_{1},\mu_{2})^2+|z_{1}-z_{2}|^{2}\).
\de
\end{enumerate}
\begin{enumerate}[$(\mathbf{H}^2_{b_{2}, \s_{2}})$]
\item
There exists a constant $\b>0$ with $\b>5L_{b_{2}, \s_{2}}$ such that for $\mu\in\cP_{2}(\mR^n)$, $z_{i}\in\mR^m$, $i=1, 2$,
\ce
&&2\<z_{1}-z_{2},b_{2}(\mu, z_{1})-b_{2}(\mu, z_{2})\>+3\|\s_{2}(\mu, z_{1})-\s_{2}(\mu, z_{2})\|^{2}
\leq -\b|z_{1}-z_{2}|^{2}.
\de
\end{enumerate}
\begin{enumerate}[$(\mathbf{H}^{2'}_{b_{2}, \s_{2}})$]
\item
For some $p>2$, there exists a constant $\b'>0$ with $\b'>(2p+1)L_{b_{2}, \s_{2}}$ such that for $\mu\in\cP_{2}(\mR^n)$, $z_{i}\in\mR^m$, $i=1, 2$,
\ce
&&2\<z_{1}-z_{2},b_{2}(\mu,z_{1})-b_{2}(\mu,z_{2})\>+(2p-1)\|\s_{2}(\mu,z_{1})-\s_{2}(\mu,z_{2})\|^{2}
\leq -\b'|z_{1}-z_{2}|^{2}.
\de
\end{enumerate}
\begin{enumerate}[$(\mathbf{H}^3_{b_{2}, \s_{2}})$]
\item
$\partial_{\mu}b_2(\mu,z)(y), \partial_z b_2(\mu,z), \partial_{zz} b_2(\mu,z), \partial_y\partial_{\mu}b_2(\mu,z)(y)$ exist for any $(\mu,y,z)\in\cP_{2}(\mR^n)\times\mR^n\times\mR^m$ and $\partial_{zz} b_2(\mu,z), \partial_y\partial_{\mu}b_2(\mu,z)(y)$ are uniformly bounded. And there exists a constant $\g_2\in(0,1]$ such that for $z_{i}\in\mR^m$, $i=1, 2$,
\ce
&&\sup\limits_{\mu\in\cP_{2}(\mR^n)}\|\partial_{\mu}b_2(\mu,z_1)-\partial_{\mu}b_2(\mu,z_2)\|_{L^2(\mu)}\leq C|z_1-z_2|^{\g_2},\\
&&\sup\limits_{\mu\in\cP_{2}(\mR^n)}\|\partial_{z}b_2(\mu,z_1)-\partial_{z}b_2(\mu,z_2)\|\leq C|z_1-z_2|^{\g_2},\\
&&\sup\limits_{\mu\in\cP_{2}(\mR^n)}\|\partial_{zz} b_2(\mu,z_1)-\partial_{zz} b_2(\mu,z_2)\|\leq C|z_1-z_2|^{\g_2},\\
&&\sup\limits_{\mu\in\cP_{2}(\mR^n)}\|\partial_y\partial_{\mu}b_2(\mu,z_1)-\partial_y\partial_{\mu}b_2(\mu,z_2)\|_{L^2(\mu)}\leq C|z_1-z_2|^{\g_2}.
\de
Moreover, the above conditions also hold for $\s_2$. 
\end{enumerate}
\begin{enumerate}[$(\mathbf{H}_{h})$]
\item $h$ is bounded, and there is a constant $L_{h}>0$ such that
$$
|h(x_{1},\mu_{1},z_{1})-h(x_{2},\mu_{2},z_{2})|^{2}\leq L_{h}(|x_{1}-x_{2}|^{2}+\mW_{2}(\mu_{1},\mu_{2})^2+|z_{1}-z_{2}|^{2}).
$$
\end{enumerate}

\br
(i) $(\mathbf{H}^1_{b_{1}, \s_{1}})$ yields that there exists a constant $\bar{L}_{b_{1}, \s_{1}}>0$ such that for $x\in\mR^n$, $\mu\in\cP_{2}(\mR^n)$, $z\in\mR^m$,
\be
|b_{1}(x,\mu,z)|^{2}+\|\s_{1}(x,\mu,z)\|^{2}
\leq \bar{L}_{b_{1}, \s_{1}}(1+|x|^{2}+\|\mu\|^{2}+|z|^{2}).
\label{b1line}
\ee

(ii) $(\mathbf{H}^1_{b_{2}, \s_{2}})$ implies that there exists a constant $\bar{L}_{b_{2}, \s_{2}}>0$ such that for $\mu\in\cP_{2}(\mR^n)$, $z\in\mR^m$,
\be
|b_{2}(\mu,z)|^{2}+\|\s_{2}(\mu,z)\|^{2}
\leq \bar{L}_{b_{2}, \s_{2}}(1+\|\mu\|^{2}+|z|^{2}).
\label{b2nu}
\ee

(iii) $(\mathbf{H}^1_{b_{2}, \s_{2}})$ and $(\mathbf{H}^2_{b_{2}, \s_{2}})$ yield that there exists a constant $C>0$ such that for $\mu\in\cP_{2}(\mR^n)$, $z\in\mR^m$
\be
&&2\<z,b_{2}(\mu,z)\>+3\|\s_{2}(\mu,z)\|^{2}\leq -\a|z|^{2}+C(1+\|\mu\|^{2}),
\label{bemu}
\ee
where $\a:=\b-4L_{b_{2}, \s_{2}}>L_{b_{2}, \s_{2}}$.

(iv) $(\mathbf{H}^1_{b_{2}, \s_{2}})$ and $(\mathbf{H}^{2'}_{b_{2}, \s_{2}})$ imply that there exists a constant $C>0$ such that for $\mu\in\cP_{2}(\mR^n)$, $z\in\mR^m$
\be
&&2\<z,b_{2}(\mu,z)\>+(2p-1)\|\s_{2}(\mu,z)\|^{2}\leq -\a'|z|^{2}+C(1+\|\mu\|^{2}),
\label{bemup}
\ee
where $\a':=\b'-2pL_{b_{2}, \s_{2}}>L_{b_{2}, \s_{2}}$.

Here we mention that $(\mathbf{H}^{2'}_{b_{2}, \s_{2}})$ is stronger than $(\mathbf{H}^2_{b_{2}, \s_{2}})$, and $(\mathbf{H}^2_{b_{2}, \s_{2}})$ and $(\mathbf{H}^{2'}_{b_{2}, \s_{2}})$  are used in Theorem \ref{weakconv} and \ref{convfilt}, respectively.
\er

\section{Main results}\label{mainresu}

In this section, we sate our main results in this paper.

\subsection{The average principle}\label{fram}

Under $(\mathbf{H}^1_{b_{1}, \s_{1}})$ $(\mathbf{H}^1_{b_{2}, \s_{2}})$, by \cite[Theorem 2.1]{WangFY} or \cite[Theorem 3.1]{DQ1}, we know that the slow-fast system (\ref{slfa}) has a unique strong solution $(X_{\cdot}^{\e}, Z_{\cdot}^{\e})$. 

Next, we take any $\mu\in\cP_{2}(\mR^n)$, and fix it. Consider the following SDE:
\be\left\{\begin{array}{l}
\dif Z_{t}^{\mu,z_0}=b_{2}(\mu,Z_{t}^{\mu,z_0})\dif t+\s_{2}(\mu,Z_{t}^{\mu,z_0})\dif W_{t},\\
Z_{0}^{\mu,z_0}=z_0, \quad 0 \leq t \leq T.
\end{array}
\right.
\label{Eq2}
\ee
Based on \cite{hzy}, it holds that under $(\mathbf{H}^1_{b_{2}, \s_{2}})$, the above equation has a unique strong solution $Z_{\cdot}^{\mu,z_0}$. Besides, $Z_{\cdot}^{\mu,z_0}$ is a Markov process. Let $p_{t}(\mu;z_0,A)$ denote the transition probability of $Z_{\cdot}^{\mu,z_0}$ for $t\geq0$ and $A\in\sB(\mR^{m})$. The associated transition semigroup $\{P_{t}^{\mu}\}_{t\geq0}$ is given by
$$
(P_{t}^{\mu}\phi)(z_0)=\int_{\mR^{m}}\phi(z')p_{t}(\mu;z_0,\dif z'), \quad \phi\in\cB_b(\mR^{m}).
$$
Under $(\mathbf{H}^2_{b_{2}, \s_{2}})$, by \cite[Theorem 3.1]{WangFY}, one could obtain that there exists a unique invariant probability measure $\nu^{\mu}$ for $Z^{\mu,z_0}_{\cdot}$. 

In the following, for $(x,\mu)\in\mR^n\times\cP_{2}(\mR^n)$, set 
\be
\bar{b}_{1}(x,\mu):=\int_{\mR^{m}}b_{1}(x,\mu,z)\nu^{\mu}(\dif z), \quad \Sigma(x,\mu):=\int_{\mR^{m}}(\s_{1}\s_{1}^{*})(x,\mu,z)\nu^{\mu}(\dif z),
\label{bbar}
\ee
and $\Sigma(x,\mu)$ is a positive definite symmetric matrix. Thus, we construct an average equation on $(\Omega,\sF,\{\sF_{t}\}_{t\in[0,T]},\mP)$ as follows:
\be\left\{\begin{array}{l}
\dif \bar{X_{t}}=\bar{b}_{1}(\bar{X_{t}},\sL^{\mP}_{\bar{X_{t}}})\dif t+\bar{\s}_1(\bar{X_{t}},\sL^{\mP}_{\bar{X_{t}}})\dif B_{t},\\
\bar{X_{0}}=x_0,
\end{array}
\right.
\label{Eq3}
\ee
where $\bar{\s}_1(x,\mu)$ is the square root of $\Sigma(x,\mu)$, i.e. $\bar{\s}_1(x,\mu)$ is a positive definite symmetric matrix  satisfying $\Sigma(x,\mu)=(\bar{\s}_1\bar{\s}_1)(x,\mu)$. 

Now, it is the position to state the first main result of this paper.

\bt\label{weakconv}
Suppose that $(\mathbf{H}^1_{b_{1}, \s_{1}})$ $(\mathbf{H}^2_{\s_{1}})$ $(\mathbf{H}^3_{b_{1}, \s_{1}})$ $(\mathbf{H}^1_{b_{2}, \s_{2}})$-$(\mathbf{H}^3_{b_{2}, \s_{2}})$ hold. Then $\{X_{t}^{\e}, t\in[0,T]\}$ converges weakly to $\{\bar{X}_{t}, t\in[0,T]\}$ in $C([0,T],\mR^{n})$, where $\bar{X_{\cdot}}$ is a strong solution of Eq.(\ref{Eq3}).
\et

The proof of the above theorem is placed in Section \ref{prooaver}.

\subsection{Convergence for nonlinear filtering}\label{Nf}

Set
\ce
(\Lambda_{t}^{\e})^{-1}:=\exp\left\{-\int_{0}^{t}h^{i}(X_{s}^{\e},\sL^{\mP}_{X_{s}^{\e}},Z_{s}^{\e})\dif V_{s}^{i} -\frac{1}{2}\int_{0}^{t}|h(X_{s}^{\e},\sL^{\mP}_{X_{s}^{\e}},Z_{s}^{\e})|^{2}\dif s\right\}.
\de
Here and hereafter, we use the convention that repeated indices imply the summation. Under $(\mathbf{H}_{h})$, we get that
$$
\mE\left(\exp\left\{\frac{1}{2}\int_{0}^{T}|h(X_{s}^{\e},\sL^{\mP}_{X_{s}^{\e}},Z_{s}^{\e})|^{2}\dif s\right\}\right)<\infty,
$$
and furthermore $(\Lambda_{t}^{\e})^{-1}$ is an exponential martingale under the measure $\mP$. Define a probability measure $\mP^{\e}$ via
\ce
\frac{\dif\mP^{\e}}{\dif \mP}=(\Lambda_{T}^{\e})^{-1}.
\de
Then by the Girsanov theorem, it holds that $Y_{\cdot}^{\e}$ is a Brownian motion under the probability measure $\mP^{\e}$.

Define the nonlinear filtering for $(X_{t}^{\e}, \sL^{\mP}_{X_{t}^{\e}}, Z_{t}^{\e})$: for any $\Psi\in \cB_{b}(\mR^{n}\times\cP_{2}(\mR^{n})\times\mR^{m})$
\ce
\rho_{t}^{\e}(\Psi)
:=\mE^{\mP^{\e}}[\Psi(X_{t}^{\e},\sL^{\mP}_{X_{t}^{\e}},Z_{t}^{\e})\Lambda_{t}^{\e}|\mathscr{F}_{t}^{Y^{\e}}], \quad \pi_{t}^{\e}(\Psi)
:=\mE[\Psi(X_{t}^{\e},\sL^{\mP}_{X_{t}^{\e}},Z_{t}^{\e})|\mathscr{F}_{t}^{Y^{\e}}],
\de
where $\mathscr{F}_{t}^{Y^{\e}}:=\sigma\{Y_{s}^{\e},0\leq s \leq t\}\vee \cN$, and $\cN$ denotes  the collection of all zero sets under the $\mP$-measure. Here $\rho_{t}^{\e}$, $\pi_{t}^{\e}$ are called the unnormalized and normalized filtering of $(X_{t}^{\e}, \sL^{\mP}_{X_{t}^{\e}}, Z_{t}^{\e})$  with respect to $\mathscr{F}_{t}^{Y^{\e}}$, respectively. By the Kallianpur-Striebel formula, we get the following relationship between $\rho_{t}^{\e}(\Psi)$ and $\pi_{t}^{\e}(\Psi)$:
$$
\pi_{t}^{\e}(\Psi)=\frac{\rho_{t}^{\e}(\Psi)}{\rho_{t}^{\e}(1)}.
$$
Define the $(x,\mu)$-marginal of $\rho_{t}^{\e}$ as follows:
\ce
\rho_{t}^{\e,x,\mu}(F)=\int_{\mR^{n}\times\cP_{2}(\mR^{n})\times\mR^{m}}F(x,\mu)\rho_{t}^{\e}(\dif x,\dif \mu,\dif z), \quad F\in \cB_{b}(\mR^{n}\times\cP_{2}(\mR^{n})),
\de
and it holds that
$$
\pi_{t}^{\e,x,\mu}(F)=\frac{\rho_{t}^{\e,x,\mu}(F)}{\rho_{t}^{\e,x,\mu}(1)},
$$
where $\pi_{t}^{\e,x,\mu}$ is the $(x,\mu)$-marginal of $\pi_{t}^{\e}$.

Next, set
\ce
&&\bar{\Lambda}_{t}:=\exp\left\{\int_{0}^{t}\bar{h}^{i}(\bar{X_{s}},\sL^{\mP}_{\bar{X_{s}}})\dif Y_{s}^{\e, i} -\frac{1}{2}\int_{0}^{t}|\bar{h}(\bar{X_{s}},\sL^{\mP}_{\bar{X_{s}}})|^{2}\dif s\right\},\\
&&\bar{\rho}_{t}(F):=\mE^{\mP^{\e}}[F(\bar{X_{t}},\sL^{\mP}_{\bar{X_{t}}})\bar{\Lambda}_{t}|\mathscr{F}_{t}^{Y^{\e}}],\\
&&\bar{\pi}_{t}(F):=\frac{\bar{\rho}_{t}(F)}{\bar{\rho}_{t}(1)},
\de
where $\bar{h}(x,\mu):=\int_{\mR^{m}}h(x,\mu,z)\nu^{\mu}(\dif z)$. Then about the relationship between $\pi_{t}^{\e,x,\mu}$ and $\bar{\pi}_{t}$, we have the following result which is the second main result for this paper.

\bt\label{convfilt}
Assume that $(\mathbf{H}^1_{b_{1}, \s_{1}})$ $(\mathbf{H}^2_{\s_{1}})$ $(\mathbf{H}^3_{b_{1}, \s_{1}})$ $(\mathbf{H}^{1}_{b_{2}, \s_{2}})$ $(\mathbf{H}^{2'}_{b_{2}, \s_{2}})$ $(p\geq 12)$ $(\mathbf{H}^{3}_{b_{2}, \s_{2}})$ and $(\mathbf{H}_{h})$ hold. Then for any $t\in[0,T]$, $\pi_{t}^{\e,x,\mu}$ converges weakly to $\bar{\pi}_{t}$ as $\e\rightarrow 0$.
\et

We will prove the above theorem in Section \ref{proofilt}.

\section{Proof of Theorem \ref{weakconv}}\label{prooaver}

In this section, we prove Theorem \ref{weakconv}. We divide the proof into four parts. In the first part (Subsection \ref{rwcx}) we show relatively weak compactness for $\{X_{t}^{\e}, t\in[0,T]\}$. In the second (Subsection \ref{soesaveq}) and third parts (Subsection \ref{soespoeq}), we present some estimates for the average equation (\ref{Eq3}) and a Poisson equation. Finally, we prove that the weak limit of $\{X_{t}^{\e}, t\in[0,T]\}$ is the unique solution $\bar{X}$ of Eq.(\ref{Eq3}) in the fourth part (Subsection \ref{weconxe}). 

\subsection{Relatively weak compactness for $\{X_{t}^{\e}, t\in[0,T]\}$}\label{rwcx}

\bl \label{xezeb}
 Under $(\mathbf{H}^{1}_{b_{1}, \s_{1}})$ $(\mathbf{H}^{1}_{b_{2}, \s_{2}})$-$(\mathbf{H}^{2}_{b_{2}, \s_{2}})$, there exists a constant $C>0$ such that 
\ce
\sup\limits_{t\in[0,T]}\mE|X_{t}^{\e}|^{4}\leq C(1+|x_0|^4+|z_0|^4), \quad \sup\limits_{t\in[0,T]}\mE|Z_{t}^{\e}|^{4}\leq C(1+|x_0|^4+|z_0|^4).
\de
\el
\begin{proof}
For $X_{t}^{\e}$, based on the H\"older inequality, the isometric formula and $(\ref{b1line})$, we can get
\be
\mE|X_{t}^{\e}|^{4}
&\leq& 3^3|x_0|^{4}+3^3\mE\Big|\int_{0}^{t}b_{1}(X_{s}^{\e},\sL^{\mP}_{X_{s}^{\e}},Z_{s}^{\e})\dif s\Big|^{4}
+3^3\mE\Big|\int_{0}^{t}\s_{1}(X_{s}^{\e},\sL^{\mP}_{X_{s}^{\e}},Z_{s}^{\e})\dif B_{s}\Big|^{4}\no\\
&\leq& 3^3|x_0|^{4}+3^3t^3\int_{0}^{t}\mE|b_{1}(X_{s}^{\e},\sL^{\mP}_{X_{s}^{\e}},Z_{s}^{\e})|^4\dif s+3^3t\int_{0}^{t}\mE\|\s_{1}(X_{s}^{\e},\sL^{\mP}_{X_{s}^{\e}},Z_{s}^{\e})\|^{4}\dif s\no\\
&\leq& 3^3|x_0|^{4}+3^3(T^3+T)\int_{0}^{t}C(1+\mE|X_{s}^{\e}|^4+\|\sL^{\mP}_{X_{s}^{\e}}\|^4+\mE|Z_{s}^{\e}|^4)\dif s\no\\
&\leq& 3^3|x_0|^{4}+3^3(T^3+T)\int_{0}^{t}C(1+2\mE|X_{s}^{\e}|^4+\mE|Z_{s}^{\e}|^4)\dif s,
\label{exqc}
\ee
where $\|\sL^{\mP}_{X_{s}^{\e}}\|^2=\mE|X_{s}^{\e}|^2$.

For $Z_{t}^{\e}$, applying the It\^{o} formula to $|Z_{t}^{\e}|^{4}$ and taking the expectation, one could obtain that
\ce
\mE|Z_{t}^{\e}|^{4}
&=&|z_0|^{4}+\frac{4}{\e}\mE\int_{0}^{t}|Z_{s}^{\e}|^2\<Z_{s}^{\e}, b_{2}(\sL^{\mP}_{X_{s}^{\e}},Z_{s}^{\e})\>\dif s\\
&&+\frac{4}{\e}\mE\int_{0}^{t}|\s_{2}(\sL^{\mP}_{X_{s}^{\e}},Z_{s}^{\e})Z_{s}^{\e}|^2\dif s+\frac{2}{\e}\mE\int_{0}^{t}|Z_{s}^{\e}|^2\|\s_{2}(\sL^{\mP}_{X_{s}^{\e}},Z_{s}^{\e})\|^{2}\dif s,
\de
and 
\ce
\frac{\dif}{\dif t}\mE|Z_{t}^{\e}|^{4}&\leq&\frac{4}{\e}\mE|Z_{t}^{\e}|^2\<Z_{t}^{\e}, b_{2}(\sL^{\mP}_{X_{t}^{\e}},Z_{t}^{\e})\>+\frac{6}{\e}\mE|Z_{t}^{\e}|^2\|\s_{2}(\sL^{\mP}_{X_{t}^{\e}},Z_{t}^{\e})\|^{2}\\
&\leq&\frac{2}{\e}\mE|Z_{t}^{\e}|^2\(-\a|Z_{t}^{\e}|^{2}+C(1+\|\sL^{\mP}_{X_{t}^{\e}}\|^{2})\)\\
&\leq&\frac{2}{\e}\left[-\a\mE|Z_{t}^{\e}|^{4}+L_{b_2,\s_2}\mE|Z_{t}^{\e}|^{4}+C(1+\|\sL^{\mP}_{X_{t}^{\e}}\|^{4})\right]\\
&\leq&\frac{-2(\a-L_{b_2,\s_2})}{\e}\mE|Z_{t}^{\e}|^{2}+\frac{C}{\e}\left(1+\mE|X_{t}^{\e}|^{4}\right),
\de
where the second and fourth inequalities are based on (\ref{bemu}) and $\|\sL^{\mP}_{X_{t}^{\e}}\|^{2}=\mE|X_{t}^{\e}|^{2}$, respectively. By the comparison theorem, we have that
\be
\mE|Z_{t}^{\e}|^{2}&\leq& |z_0|^{2}e^{-\frac{2(\a-L_{b_2,\s_2})}{\e}t}+\frac{C}{\e}\int_{0}^{t}e^{-\frac{2(\a-L_{b_2,\s_2})}{\e}(t-s)}\left(1+\mE|X_{s}^{\e}|^{4}\right)\dif s\no\\
&\leq&|z_0|^{4}+C\left(1+\sup\limits_{s\in[0,t]}\mE|X_{s}^{\e}|^{4}\right).
\label{zees}
\ee

Finally, inserting (\ref{zees}) in (\ref{exqc}), by the Gronwall inequality one can get that
\ce
\sup\limits_{t\in[0,T]}\mE|X_{t}^{\e}|^{4}\leq C(1+|x_0|^4+|z_0|^4),
\de
which is just the required first result. Moreover, this together with (\ref{zees}) implies the second estimate.
\end{proof}

\bl \label{xediff}
 Suppose that $(\mathbf{H}^{1}_{b_{1}, \s_{1}})$ $(\mathbf{H}^{1}_{b_{2}, \s_{2}})$-$(\mathbf{H}^{2}_{b_{2}, \s_{2}})$ hold. Then it holds that for any $\delta>0$ and any $\{\mathscr{F}_{t}\}$-stopping time $\tau\leq T$,
 \ce
 \mE|X_{\tau+\d}^{\e}-X_{\tau}^{\e}|^{2}\leq C(1+|x_0|^2+|z_0|^2)(\d^{2}+\d).
 \de
 \el
 \begin{proof}
From the H\"older inequality, the isometric formula and Lemma \ref{xezeb}, it follows that
\ce
\mE|X_{\tau+\d}^{\e}-X_{\tau}^{\e}|^{2}
&\leq& 2\mE\Big|\int_{\tau}^{\tau+\d}b_{1}(X_{s}^{\e},\sL^{\mP}_{X_{s}^{\e}},Z_{s}^{\e})\dif s\Big|^{2}
+2\mE\Big|\int_{\tau}^{\tau+\d}\s_{1}(X_{s}^{\e},\sL^{\mP}_{X_{s}^{\e}},Z_{s}^{\e})\dif B_{s}\Big|^{2}\\
&\leq& 2\d\mE\int_{\tau}^{\tau+\d}|b_{1}(X_{s}^{\e},\sL^{\mP}_{X_{s}^{\e}},Z_{s}^{\e})|^{2}\dif s+2\mE\int_{\tau}^{\tau+\d}\|\s_{1}(X_{s}^{\e},\sL^{\mP}_{X_{s}^{\e}},Z_{s}^{\e})\|^{2}\dif s\\
&\leq& (2\d+2)C\mE\int_{\tau}^{\tau+\d}(1+|X_{s}^{\e}|^2+\|\sL^{\mP}_{X_{s}^{\e}}\|^2+|Z_{s}^{\e}|^{2})\dif s\\
&=& (2\d+2)C\mE\int_{\tau}^{\tau+\d}(1+|X_{s}^{\e}|^2+\mE|X_{s}^{\e}|^2+|Z_{s}^{\e}|^{2})\dif s\\
&\leq&(2\d+2)C\mE\int_{\tau}^{\tau+\d}(1+|X_{s}^{\e}|^2+C(1+|x_0|^2+|z_0|^2)+|Z_{s}^{\e}|^{2})\dif s\\
&=&(2\d+2)C\mE\int_{0}^{\d}(1+|X_{\tau+u}^{\e}|^2+C(1+|x_0|^2+|z_0|^2)+|Z_{\tau+u}^{\e}|^{2})\dif u\\
&\leq& C(1+|x_0|^2+|z_0|^2)(\d^{2}+\d),
\de
which completes the proof.
\end{proof}

Now, combining Lemma \ref{xezeb} and \ref{xediff} with \cite[Theorem 2.7]{TGK}, we draw the following conclusion.

\bp \label{xeweakcomp}
 Under $(\mathbf{H}^{1}_{b_{1}, \s_{1}})$, $(\mathbf{H}^{1}_{b_{2}, \s_{2}})$-$(\mathbf{H}^{2}_{b_{2}, \s_{2}})$, $\{X_{t}^{\e}, t\in[0,T]\}$ is relatively weakly compact in $C([0,T],\mR^{n})$.
\ep

\subsection{Some estimates for the average equation (\ref{Eq3})}\label{soesaveq}

Since the average equation (\ref{Eq3}) depends on the invariant probability measure of the frozen equation (\ref{Eq2}), we begin with a result about Eq.(\ref{Eq2}).

\bl\label{frozesti}
Assume that $(\mathbf{H}^{1}_{b_{2}, \s_{2}})$-$(\mathbf{H}^{2}_{b_{2}, \s_{2}})$ hold. Then we have for $\mu, \mu_1, \mu_2\in\cP_2(\mR^n), z, z_1, z_2\in\mR^m$,
\ce
&&\mE|Z_{t}^{\mu,z}|^2\leq |z|^2e^{-\a t}+C(1+\|\mu\|^2), \quad t\geq 0,\\
&&\mE|Z_{t}^{\mu_{1},z_1}-Z_{t}^{\mu_{2},z_2}|^{2}\leq|z_1-z_2|^2e^{-\a t}+\frac{C}{\a}\mW_2(\mu_1,\mu_2)^2, \quad t\geq 0.
\de
\el
\begin{proof}
First of all, by the same deduction to that for $Z^\e$ in Lemma \ref{xezeb}, one can get the first estimate. 

For the second estimate, take $T>0$ and note that for any $t\in[0,T]$
\ce
Z_{t}^{\mu_{1},z_1}-Z_{t}^{\mu_{2},z_2}
&=&\int_{0}^{t}\(b_{2}(\mu_{1},Z_{s}^{\mu_{1},z_1})-b_{2}(\mu_{2},Z_{s}^{\mu_{2},z_2})\)\dif s\\
&&+\int_{0}^{t}\(\s_{2}(\mu_{1},Z_{s}^{\mu_{1},z_1})-\s_{2}(\mu_{2},Z_{s}^{\mu_{2},z_2})\)\dif W_{s}.
\de
Thus, applying the It\^{o} formula to $|Z_{t}^{\mu_{1},z_1}-Z_{t}^{\mu_{2},z_2}|^{2}$ and taking expectation on two sides, we obtain that
\ce
&&\mE|Z_{t}^{\mu_{1},z_1}-Z_{t}^{\mu_{2},z_2}|^{2}\\
&=&2\mE\int_{0}^{t}\< Z_{s}^{\mu_{1},z_1}-Z_{s}^{\mu_{2},z_2}, b_{2}(\mu_{1},Z_{s}^{\mu_{1},z_1})-b_{2}(\mu_{2},Z_{s}^{\mu_{2},z_2})\>\dif s\\
&&+\mE\int_{0}^{t}\|\s_{2}(\mu_{1},Z_{s}^{\mu_{1},z_1})-\s_{2}(\mu_{2},Z_{s}^{\mu_{2},z_2})\|^{2}\dif s,
\de
and 
\ce
&&\frac{\dif \mE|Z_{t}^{\mu_{1},z_1}-Z_{t}^{\mu_{2},z_2}|^{2}}{\dif t}\\ 
&\leq&2\mE\< Z_{t}^{\mu_{1},z_1}-Z_{t}^{\mu_{2},z_2}, b_{2}(\mu_{1},Z_{t}^{\mu_{1},z_1})-b_{2}(\mu_{1},Z_{t}^{\mu_{2},z_2})\>+\mE\|\s_{2}(\mu_{1},Z_{t}^{\mu_{1},z_1})-\s_{2}(\mu_{1},Z_{t}^{\mu_{2},z_2})\|^{2}\\
&&+2\mE\< Z_{t}^{\mu_{1},z_1}-Z_{t}^{\mu_{2},z_2}, b_{2}(\mu_{1},Z_{t}^{\mu_{2},z_2})-b_{2}(\mu_{2},Z_{t}^{\mu_{2},z_2})\>\\
&&+\mE\|\s_{2}(\mu_{1},Z_{t}^{\mu_{1},z_1})-\s_{2}(\mu_{1},Z_{t}^{\mu_{2},z_2})\|^{2}+2\mE\|\s_{2}(\mu_{1},Z_{t}^{\mu_{2},z_2})-\s_{2}(\mu_{2},Z_{t}^{\mu_{2},z_2})\|^{2}\\
&\leq&-\b\mE| Z_{t}^{\mu_{1},z_1}-Z_{t}^{\mu_{2},z_2}|^2\dif t+2\mE|Z_{t}^{\mu_{1},z_1}-Z_{t}^{\mu_{2},z_2}||b_{2}(\mu_{1},Z_{t}^{\mu_{2},z_2})-b_{2}(\mu_{2},Z_{t}^{\mu_{2},z_2})|\\
&&+L_{b_{2}, \s_{2}}\mE| Z_{t}^{\mu_{1},z_1}-Z_{t}^{\mu_{2},z_2}|^2+2L_{b_{2}, \s_{2}}\mW_2(\mu_1,\mu_2)^2\\
&\leq&-\a\mE| Z_{t}^{\mu_{1},z_1}-Z_{t}^{\mu_{2},z_2}|^2+C\mW_2(\mu_1,\mu_2)^2.
\de
By the comparison theorem, it holds that
$$
\mE|Z_{t}^{\mu_{1},z_1}-Z_{t}^{\mu_{2},z_2}|^{2}\leq|z_1-z_2|^2e^{-\a t}+\frac{C}{\a}\mW_2(\mu_1,\mu_2)^2,
$$
which completes the proof.
\end{proof}

By the above lemma and the definition of $\nu^{\mu}$, it holds that
\ce
&&\int_{\mR^{m}}|z|^{2}\nu^{\mu}(\dif z)=\int_{\mR^{m}}\mE|Z_{t}^{\mu,z}|^{2}\nu^{\mu}(\dif z)\leq \int_{\mR^{m}}\(|z|^{2}e^{-\a t}+C(1+\|\mu\|^{2})\)\nu^{\mu}(\dif z)\no\\
&=& e^{-\a t}\int_{\mR^{m}}|z|^{2}\nu^{\mu}(\dif z)+C(1+\|\mu\|^{2}),
\de
and furthermore
\be
\int_{\mR^{m}}|z|^{2}\nu^{\mu}(\dif z)\leq C(1+\|\mu\|^{2}).
\label{inu2}
\ee

\bl\label{barsig}
 Suppose that $(\mathbf{H}^1_{b_{1}, \s_{1}})$, $(\mathbf{H}^2_{\s_{1}})$ and $(\mathbf{H}^{1}_{b_{2}, \s_{2}})$-$(\mathbf{H}^2_{b_{2}, \s_{2}})$  hold. Then there exists a constant $C>0$ such that
\ce
\|\bar{\s}_1(x_{1},\mu_{1})-\bar{\s}_1(x_{2},\mu_{2})\|\leq C(|x_{1}-x_{2}|+\mW_2(\mu_{1},\mu_{2})).
\de
Moreover, Eq.(\ref{Eq3}) has a unique strong solution $\bar{X}$ satisfying for any $q\geq 1$
\ce
\sup\limits_{t\in[0,T]}\mE|\bar{X}_{t}|^{2q}\leq C(1+|x_0|^{2q}).
\de
\el
\begin{proof}
First of all, note that
\be
\|\bar{\s}_1(x_{1},\mu_{1})-\bar{\s}_1(x_{2},\mu_{2})\|
&=&\|(\Sigma(x_{1},\mu_{1})-\Sigma(x_{2},\mu_{2}))(\bar{\s}_1(x_{1},\mu_{1})+\bar{\s}_1(x_{2},\mu_{2}))^{-1}\|\no\\
&\leq& \|\Sigma(x_{1},\mu_{1})-\Sigma(x_{2},\mu_{2})\|\|(\bar{\s}_1(x_{1},\mu_{1})+\bar{\s}_1(x_{2},\mu_{2}))^{-1}\|\no\\
&\leq& C\|\Sigma(x_{1},\mu_{1})-\Sigma(x_{2},\mu_{2})\|,
\label{sSig}
\ee
where the last inequality is based on $(\mathbf{H}^2_{\s_{1}})$. Therefore, we only need to estimate $\|\Sigma(x_{1},\mu_{1})-\Sigma(x_{2},\mu_{2})\|$. From the definition of $\Sigma(x,\mu)$, the H\"{o}lder inequality  and Lemma \ref{frozesti}, it follows that
\ce
&&\|\Sigma(x_{1},\mu_{1})-\Sigma(x_{2},\mu_{2})\|\no\\
&=&\left\|\int_{\mR^{m}}(\s_{1}\s_{1}^{*})(x_{1},\mu_{1},z)\nu^{\mu_{1}}(\dif z)
-\int_{\mR^{m}}(\s_{1}\s_{1}^{*})(x_{2},\mu_{2},z)\nu^{\mu_{2}}(\dif z)\right\|\no\\
&=&\left\|\lim_{S\rightarrow\infty}\frac{1}{S}\int_{0}^{S}\mE(\s_{1}\s_{1}^{*})(x_{1},\mu_{1},Z_{t}^{\mu_{1},z_0})\dif t
-\lim_{S\rightarrow\infty}\frac{1}{S}\int_{0}^{S}\mE(\s_{1}\s_{1}^{*})(x_{2},\mu_{2},Z_{t}^{\mu_{2},z_0})\dif t\right\|\no\\
&\leq&\lim_{S\rightarrow\infty}\frac{1}{S}\int_{0}^{S}\mE\|(\s_{1}\s_{1}^{*})(x_{1},\mu_{1},Z_{t}^{\mu_{1},z_0})
-(\s_{1}\s_{1}^{*})(x_{2},\mu_{2},Z_{t}^{\mu_{2},z_0})\|\dif t\no\\
&\leq&\lim_{S\rightarrow\infty}\frac{1}{S}\int_{0}^{S}\mE\|\s_{1}(x_{1},\mu_{1},Z_{t}^{\mu_{1},z_0})\s_{1}^{*}(x_{1},\mu_{1},Z_{t}^{\mu_{1},z_0})\no\\
&&\qquad\qquad -\s_{1}(x_{1},\mu_{1},Z_{t}^{\mu_{1},z_0})\s_{1}^{*}(x_{2},\mu_{2},Z_{t}^{\mu_{2},z_0})\|\dif t\no\\
&&+\lim_{S\rightarrow\infty}\frac{1}{S}\int_{0}^{S}\mE\|\s_{1}(x_{1},\mu_{1},Z_{t}^{\mu_{1},z_0})\s_{1}^{*}(x_{2},\mu_{2},Z_{t}^{\mu_{2},z_0})\no\\
&&\qquad\qquad -\s_{1}(x_{2},\mu_{2},Z_{t}^{\mu_{2},z_0})\s_{1}^{*}(x_{2},\mu_{2},Z_{t}^{\mu_{2},z_0})\|\dif t\no\\
&\leq&\lim_{S\rightarrow\infty}\frac{1}{S}\int_{0}^{S}\left(\mE\|\s_{1}^{*}(x_{1},\mu_{1},Z_{t}^{\mu_{1},z_0})-\s_{1}^{*}(x_{2},\mu_{2},Z_{t}^{\mu_{2},z_0})\|^2\right)^{1/2}\no\\
&&\qquad\qquad \left(\mE\|\s_{1}(x_{1},\mu_{1},Z_{t}^{\mu_{1},z_0})\|^2\right)^{1/2}\dif t\no\\
&&+\lim_{S\rightarrow\infty}\frac{1}{S}\int_{0}^{S}\left(\mE\|\s_{1}(x_{1},\mu_{1},Z_{t}^{\mu_{1},z_0})-\s_{1}(x_{2},\mu_{2},Z_{t}^{\mu_{2},z_0})\|^2\right)^{1/2}\no\\
&&\qquad\qquad \left(\mE\|\s_{1}^{*}(x_{2},\mu_{2},Z_{t}^{\mu_{2},z_0})\|^2\right)^{1/2}\dif t\no\\
&\leq& C(|x_{1}-x_{2}|+\mW_2(\mu_{1},\mu_{2}))+C\lim_{S\rightarrow\infty}\frac{1}{S}\int_{0}^{S}\left(\mE|Z_{t}^{\mu_{1},z_0}-Z_{t}^{\mu_{2},z_0}|^{2}\right)^{1/2}\dif t\\
&\leq& C(|x_{1}-x_{2}|+\mW_2(\mu_{1},\mu_{2})).
\de

Finally, inserting the above inequality in (\ref{sSig}), one could obtain that
$$
\|\bar{\s}_1(x_{1},\mu_{1})-\bar{\s}_1(x_{2},\mu_{2})\|\leq C(|x_{1}-x_{2}|+\mW_2(\mu_{1},\mu_{2})).
$$
Besides, by $(\mathbf{H}^1_{b_{1}, \s_{1}})$ and the similar deduction to that in \cite[Lemma 3.8]{SunX}, it holds that
\be
|\bar{b}_{1}(x_{1},\mu_{1})-\bar{b}_{1}(x_{2},\mu_{2})|^2
\leq L_{\bar{b}_{1}}\(|x_{1}-x_{2}|^{2}+\mW_2(\mu_{1},\mu_{2})^2\),
\label{baxm}
\ee
where $L_{\bar{b}_{1}}>0$ is a constant. Thus, Eq.(\ref{Eq3}) has a unique strong solution $\bar{X_{\cdot}}$ (c.f. \cite[Theorem 3.1]{DQ1}). Then, by the same deduction to that in Lemma \ref{xezeb}, we also obtain the required estimate. This proof is complete.
\end{proof}

\subsection{Some estimates for a Poisson equation}\label{soespoeq}

First of all, we notice that the infinitesimal generator of $Z^{\mu,z_0}_{\cdot}$ is as follows:
\ce
(\cL_0\phi)(\mu,z)
:&=&\partial_{z_{i}}\phi(z)b_{2}^{i}(\mu,z)
+\frac{1}{2}\partial_{z_{i}z_{k}}\phi(z)(\s_{2}\s_{2}^{*})^{ik}(\mu,z), \quad \phi\in C^2(\mR^m).
\de
Consider the following Poisson equation: for any $F\in C^{2,(1,1)}(\mR^n\times\cP_2(\mR^n))$
\be
(\cL_0\chi)(x,\mu,z)=-\left[(\cL-\bar{\cL})F\right](x,\mu,z), 
\label{pois}
\ee
where 
\ce
(\cL F)(x,\mu,z):&=&\partial_{x_{i}}F(x,\mu)b_{1}^{i}(x,\mu,z)+\frac{1}{2}\partial_{x_{i}x_{j}}F(x,\mu)(\s_{1}\s_{1}^{*})^{ij}(x,\mu,z)\\
&&+\int_{\mR^{n}}(\partial_{\mu}F)_{i}(x,\mu)(y)b_{1}^{i}(y,\mu,z)\mu(\dif y)\\
&&+\frac{1}{2}\int_{\mR^{n}}\partial_{y_{i}}(\partial_{\mu}F)_{j}(x,\mu)(y)(\s_{1}\s_{1}^{*})^{ij}
(y,\mu,z)\mu(\dif y),
\de
and
\ce
(\bar{\cL}F)(x,\mu)
:&=&\partial_{x_{i}}F(x,\mu)\bar{b}_{1}^{i}(x,\mu)+\frac{1}{2}\partial_{x_{i}x_{j}}F(x,\mu)(\bar{\s}_{1}\bar{\s}_{1}^{*})^{ij}(x,\mu)\\
&&+\int_{\mR^{n}}(\partial_{\mu}F)_{i}(x,\mu)(y)\bar{b}_{1}^{i}(y,\mu)\mu(\dif y)\\
&&+\frac{1}{2}\int_{\mR^{n}}\partial_{y_{i}}(\partial_{\mu}F)_{j}(x,\mu)(y)(\bar{\s}_{1}\bar{\s}_{1}^{*})^{ij}
(y,\mu)\mu(\dif y).
\de
Set $\Phi(x,\mu,z):=(\cL F)(x,\mu,z)$, and by the simple calculation, it holds that
$$
\bar{\Phi}(x,\mu):=\int_{\mR^m}\Phi(x,\mu,z)\nu^{\mu}(\dif z)=(\bar{\cL}F)(x,\mu).
$$

\bp\label{chib}
Suppose that $(\mathbf{H}^1_{b_{1}, \s_{1}})$, $(\mathbf{H}^2_{\s_{1}})$, $(\mathbf{H}^3_{b_{1}, \s_{1}})$, $(\mathbf{H}^{1}_{b_{2}, \s_{2}})$-$(\mathbf{H}^{3}_{b_{2}, \s_{2}})$ hold. For any $F(\cdot,\cdot)\in \mC^{4,(2,2)}_{b}(\mR^{n}\times\cP_{2}(\mR^{n}))$, set
\be
\chi_{F}(x,\mu, z)&:=&\int_{0}^{\infty}P_{t}^{\mu}\left[\Phi(x,\mu,\cdot)-\bar{\Phi}(x,\mu)\right](z)\dif t\no\\
&=&\int_{0}^{\infty}\(\mE \Phi(x, \mu, Z_{t}^{\mu,z})-\bar{ \Phi}(x, \mu)\)\dif t.
\label{chif}
\ee
Then $\chi_{F}(x,\mu, z)$ belongs to $C^{2,(1,1),2}(\mR^{n}\times\cP_{2}(\mR^{n})\times\mR^m)$ and is the unique solution for Eq.(\ref{pois}). Moreover, it holds that for $x\in\mR^n, \mu\in\cP_{2}(\mR^{n}), z\in\mR^m$
\ce
&&\max\{|\chi_F(x,\mu,z)|, |\p_{x}\chi_F(x,\mu,z)|,\|\p_{\mu}\chi_F(x,\mu,z)(\cdot)\|_{L^2(\mu)}, |\p_{z}\chi_F(x,\mu,z)|\}\\
&\leq& C(1+\|\mu\|+|z|),\\
&&\max\{\|\p_{xx}\chi_F(x,\mu,z)\|, \|\p_{u}\p_{\mu}\chi_F(x,\mu,z)(\cdot)\|_{L^2(\mu)}\}\\
&\leq& C(1+\|\mu\|+|z|).
\de
\ep
\begin{proof}
{\bf Step 1.} We prove that the right side of (\ref{chif}) is well-defined.

First of all, from $(\mathbf{H}^1_{b_{1}, \s_{1}})$ and $F(\cdot,\cdot)\in \mC^{4,(2,2)}_{b}(\mR^{n}\times\cP_{2}(\mR^{n}))$, it follows that $\Phi$ is Lipschitz continuous in $z$, which together with the definition of $\nu^{\mu}$ and Lemma \ref{frozesti} implies that
\be
&&\left|P_{t}^{\mu}\left[\Phi(x,\mu,\cdot)-\bar{\Phi}(x,\mu)\right](z)\right|^2=|\mE \Phi(x, \mu, Z_{t}^{\mu,z})-\bar{ \Phi}(x, \mu)|^2\no\\
&=&|\mE \Phi(x, \mu, Z_{t}^{\mu,z})-\int_{\mR^m}\Phi(x,\mu,y)\nu^{\mu}(\dif y)|^2\no\\
&=&|\mE \Phi(x, \mu, Z_{t}^{\mu,z})-\int_{\mR^m}\mE \Phi(x, \mu, Z_{t}^{\mu,y})\nu^{\mu}(\dif y)|^2\no\\
&\leq&\int_{\mR^m}\mE|\Phi(x, \mu, Z_{t}^{\mu,z})-\Phi(x, \mu, Z_{t}^{\mu,y})|^2\nu^{\mu}(\dif y)\no\\
&\leq&C\int_{\mR^m}|z-y|^2e^{-\a t}\nu^{\mu}(\dif y)\overset{(\ref{inu2})}{\leq}Ce^{-\a t}(1+|z|+\|\mu\|)^2.
\label{semiesti}
\ee
Thus, it holds that
\be 
&&\int_{0}^{\infty}\left|P_{t}^{\mu}\left[\Phi(x,\mu,\cdot)-\bar{\Phi}(x,\mu)\right](z)\right|\dif t=\int_{0}^{\infty}|\mE \Phi(x, \mu, Z_{t}^{\mu,z})-\bar{ \Phi}(x, \mu)|\dif t\no\\
&\leq&\int_{0}^{\infty}Ce^{-\frac{\a}{2} t}(1+\|\mu\|+|z|)\dif t=\frac{2C}{\a}(1+\|\mu\|+|z|).
\label{xfbo}
\ee
So the right side of (\ref{chif}) is well-defined. 

{\bf Step 2.} We show that $\chi_{F}(x,\mu, z)$ belongs to $C^{2,(1,1),2}(\mR^{n}\times\cP_{2}(\mR^{n})\times\mR^m)$ and is the unique solution for Eq.(\ref{pois}). 

First of all, we study the regularity of $\Phi$. By $(\mathbf{H}^3_{b_{1}, \s_{1}})$ and $F(\cdot,\cdot)\in \mC^{4,(2,2)}_{b}(\mR^{n}\times\cP_{2}(\mR^{n}))$, it holds that $\Phi$ belongs to $C^{2,(1,1),2}(\mR^{n}\times\cP_{2}(\mR^{n})\times\mR^m)$. We only compute $\partial_\mu\Phi(x,\mu,z)(u)$ as follows:
\ce
&&\partial_\mu\Phi(x,\mu,z)(u)\\
&=&\partial_\mu\partial_{x_{i}}F(x,\mu)(u)b_{1}^{i}(x,\mu,z)+\partial_{x_{i}}F(x,\mu)\partial_\mu b_{1}^{i}(x,\mu,z)(u)\\
&&+\frac{1}{2}\partial_\mu\partial_{x_{i}x_{j}}F(x,\mu)(u)(\s_{1}\s_{1}^{*})^{ij}(x,\mu,z)+\frac{1}{2}\partial_{x_{i}x_{j}}F(x,\mu)\partial_\mu(\s_{1}\s_{1}^{*})^{ij}(x,\mu,z)(u)\\
&&+\int_{\mR^{n}}(\partial^2_{\mu}F)_{i}(x,\mu)(y,u)b_{1}^{i}(y,\mu,z)\mu(\dif y)+\int_{\mR^{n}}(\partial_{\mu}F)_{i}(x,\mu)(y)\partial_\mu b_{1}^{i}(y,\mu,z)(u)\mu(\dif y)\\
&&+\int_{\mR^{n}}\partial_y(\partial_{\mu}F)_{i}(x,\mu)(y)b_{1}^{i}(y,\mu,z)\mu(\dif y)+\int_{\mR^{n}}(\partial_{\mu}F)_{i}(x,\mu)(y)\partial_y b_{1}^{i}(y,\mu,z)\mu(\dif y)\\
&&+\frac{1}{2}\int_{\mR^{n}}\partial_{y_{i}}(\partial^2_{\mu}F)_{j}(x,\mu)(y,u)(\s_{1}\s_{1}^{*})^{ij}(y,\mu,z)\mu(\dif y)\\
&&+\frac{1}{2}\int_{\mR^{n}}\partial_{y_{i}}(\partial_{\mu}F)_{j}(x,\mu)(y)\partial_\mu(\s_{1}\s_{1}^{*})^{ij}(y,\mu,z)(u)\mu(\dif y)\\
&&+\frac{1}{2}\int_{\mR^{n}}\partial_y\partial_{y_{i}}(\partial_{\mu}F)_{j}(x,\mu)(y)(\s_{1}\s_{1}^{*})^{ij}(y,\mu,z)\mu(\dif y)\\
&&+\frac{1}{2}\int_{\mR^{n}}\partial_{y_{i}}(\partial_{\mu}F)_{j}(x,\mu)(y)\partial_y(\s_{1}\s_{1}^{*})^{ij}(y,\mu,z)\mu(\dif y).
\de

Besides, note that $Z_{\cdot}^{\mu,z}$ satisfies Eq.(\ref{Eq2}). So, $(\mathbf{H}^{3}_{b_{2}, \s_{2}})$ assures the existence of $\partial_{\mu} Z_{t}^{\mu,z}(u),\\ \partial_u\partial_{\mu} Z_{t}^{\mu,z}(u), \partial_{z} Z_{t}^{\mu,z}, \partial_{zz} Z_{t}^{\mu,z}$.

Combining the above deduction, we have that $\chi_{F}(x,\mu, z)$ belongs to $C^{2,(1,1),2}(\mR^{n}\times\cP_{2}(\mR^{n})\times\mR^m)$. Then by acting the generator $\cL_0$ on $\chi_{F}(x,\mu,z)$, it holds that 
\ce
(\cL_0\chi_{F})(x,\mu,z)&=&\int_{0}^{\infty}(\cL_0P_{t}^{\mu})\left[\Phi(x,\mu,\cdot)-\bar{\Phi}(x,\mu)\right](z)\dif t\no\\
&=&\int_{0}^{\infty}\frac{\dif P_{t}^{\mu}\left[\Phi(x,\mu,\cdot)-\bar{\Phi}(x,\mu)\right](z)}{\dif t}\dif t\no\\
&=& \lim_{t\rightarrow\infty}P_{t}^{\mu}\left[\Phi(x,\mu,\cdot)-\bar{\Phi}(x,\mu)\right](z)-\left[\Phi(x,\mu,\cdot)-\bar{\Phi}(x,\mu)\right](z)\no\\
&=&-\left[\Phi(x,\mu,\cdot)-\bar{\Phi}(x,\mu)\right](z),
\de
which yields that $\chi_F$ is a solution for Eq.(\ref{pois}). Moreover, based on $(\mathbf{H}^{3}_{b_{2}, \s_{2}})$, we know that the solutions of Eq.(\ref{pois}) are unique up to an additive 
constant. Thus, $\chi_{F}(x,\mu, z)$ is the unique solution for Eq.(\ref{pois}). 

{\bf Step 3.} We establish the required estimates.

By (\ref{xfbo}), we conclude that 
$$
\max\{|\chi_F(x,\mu,z)|, |\p_{x}\chi_F(x,\mu,z)|, \|\p_{xx}\chi_F(x,\mu,z)\|\}\leq C(1+\|\mu\|+|z|).
$$
And from Lemma \ref{frozesti}, it follows that $|\p_{z}\chi_F(x,\mu,z)|\leq C$.

Next, we estimate $\p_{\mu}\chi_F(x,\mu,z)(u), \p_u\p_{\mu}\chi_F(x,\mu,z)(u)$. For any $s>0$, put 
\ce
\tilde{\Phi}_{s}(x,\mu,z,t):=\mE[\Phi(x,\mu,Z^{\mu,z}_t)]-\mE[\Phi(x,\mu,Z^{\mu,z}_{s+t})],
\de
and by the definition of $\nu^{\mu}$, it holds that
$$
\lim\limits_{s\rightarrow\infty}\tilde{\Phi}_{s}(x,\mu,z,t)=\mE[\Phi(x,\mu,Z^{\mu,z}_t)]-\bar{\Phi}(x,\mu).
$$
So, in order to estimate $\p_{\mu}\chi_F(x,\mu,z)(u), \p_u\p_{\mu}\chi_F(x,\mu,z)(u)$, we study $\p_\mu\tilde{\Phi}_{s}(x,\mu,z,t)(u), \\\p_u\p_\mu\tilde{\Phi}_{s}(x,\mu,z,t)(u)$.

By the similar deduction to that for (3.24), (4.6) in \cite{SunX}, there exists a $\eta>0$ such that
\ce
&&\|\p_\mu\tilde{\Phi}_{s}(x,\mu,z,t)\|_{L^2(\mu)}\leq Ce^{-\eta t}(1+\|\mu\|+|z|),\\
&&\|\p_u\p_\mu\tilde{\Phi}_{s}(x,\mu,z,t)\|_{L^2(\mu)}\leq Ce^{-\eta t}(1+\|\mu\|+|z|),
\de
which yields that
\ce
&&\|\p_{\mu}\chi_F(x,\mu,z)(\cdot)\|_{L^2(\mu)}\leq C(1+\|\mu\|+|z|),\\
&&\|\p_{u}\p_{\mu}\chi_F(x,\mu,z)(\cdot)\|_{L^2(\mu)}\leq C(1+\|\mu\|+|z|).
\de
The proof is complete.
\end{proof}

\subsection{Weak convergence of $\{X_{t}^{\e}, t\in[0,T]\}$ to $\bar{X}$}\label{weconxe}

\bp\label{xeconx}
Suppose that $(\mathbf{H}^1_{b_{1}, \s_{1}})$ $(\mathbf{H}^2_{\s_{1}})$ $(\mathbf{H}^3_{b_{1}, \s_{1}})$ $(\mathbf{H}^1_{b_{2}, \s_{2}})$-$(\mathbf{H}^3_{b_{2}, \s_{2}})$ hold. Then $\{X_{t}^{\e}, t\in[0,T]\}$ converges weakly to $\{\bar{X}_{t}, t\in[0,T]\}$ in $C([0,T],\mR^{n})$.
\ep
\begin{proof}
First of all, applying the It\^{o} formula to $F(X_{t}^{\e},\sL^{\mP}_{X_{t}^{\e}})$ for $F\in \mC_b^{4,(2,2)}(\mR^n\times\cP_2(\mR^n))$, one could obtain that for $0\leq s<t\leq T$
\ce
F(X_{t}^{\e},\sL^{\mP}_{X_{t}^{\e}})&=&F(X_{s}^{\e},\sL^{\mP}_{X_{s}^{\e}})+\int_{s}^{t}(\cL F)(X_{r}^{\e},\sL^{\mP}_{X_{r}^{\e}},Z_{r}^{\e})\dif r\\
&&+\int_{s}^{t}\frac{\partial F}{\partial x_{i}}(X_{r}^{\e},\sL^{\mP}_{X_{r}^{\e}})\s_{1}^{ij}(X_{r}^{\e},\sL^{\mP}_{X_{r}^{\e}},Z_{r}^{\e})\dif B_{r}^{j},
\de
and furthermore
\be
&&F(X_{t}^{\e},\sL^{\mP}_{X_{t}^{\e}})-F(X_{s}^{\e},\sL^{\mP}_{X_{s}^{\e}})
-\int_{s}^{t}(\bar{\cL}F)(X_{r}^{\e},\sL^{\mP}_{X_{r}^{\e}})\dif r\no\\
&=&\int_{s}^{t}(\cL F)(X_{r}^{\e}, \sL^{\mP}_{X_{r}^{\e}}, Z_{r}^{\e})\dif r-\int_{s}^{t}(\bar{\cL}F)(X_{r}^{\e},\sL^{\mP}_{X_{r}^{\e}})\dif r\no\\
&&+\int_{s}^{t}\frac{\partial F}{\partial x_{i}}(X_{r}^{\e},\sL^{\mP}_{X_{r}^{\e}})\s_{1}^{ij}(X_{r}^{\e},\sL^{\mP}_{X_{r}^{\e}},Z_{r}^{\e})\dif B_{r}^{j}.
\label{goper}
\ee
Thus, multiplying a bounded $\mathscr{F}_{s}$-measurable functional $\Gamma_{s}$ of the process $\{\bar{X}_{t}, t\in[0,T]\}$ and taking the expectation under the measure $\mP$ on both sides of (\ref{goper}), we have that
\be
&&\mE\Bigg[\Gamma_{s}(\bar{X}_{.})\(F(X_{t}^{\e},\sL^{\mP}_{X_{t}^{\e}})-F(X_{s}^{\e},\sL^{\mP}_{X_{s}^{\e}})
-\int_{s}^{t}(\bar{\cL}F)(X_{r}^{\e},\sL^{\mP}_{X_{r}^{\e}})\dif r\)\Bigg]\no\\
&=&\mE\Bigg[\Gamma_{s}(\bar{X}_{.})\int_{s}^{t}\[(\cL F)(X_{r}^{\e},\sL^{\mP}_{X_{r}^{\e}},Z_{r}^{\e})
-(\bar{\cL}F)(X_{r}^{\e},\sL^{\mP}_{X_{r}^{\e}})\]\dif r\Bigg]\no\\
&=&-\mE\Bigg[\Gamma_{s}(\bar{X}_{.})\int_{s}^{t}(\cL_0\chi_{F})(X_{r}^{\e},\sL^{\mP}_{X_{r}^{\e}},Z_{r}^{\e})\dif r\Bigg],
\label{limer}
\ee
where Proposition \ref{chib} is used in the last equality.

Next, we observe the right hand side of (\ref{limer}). By applying the It\^{o} formula to $\e\chi_{F}(X_{t}^{\e},\sL^{\mP}_{X_{t}^{\e}},Z_{t}^{\e})$, it holds that
\ce
&&\e\chi_{F}(X_{t}^{\e},\sL^{\mP}_{X_{t}^{\e}},Z_{t}^{\e})-\e\chi_{F}(X_{s}^{\e},\sL^{\mP}_{X_{s}^{\e}},Z_{s}^{\e})
-\e\int_{s}^{t}(\cL\chi_{F})(X_{r}^{\e},\sL^{\mP}_{X_{r}^{\e}},Z_{r}^{\e})\dif r\\
&=&\e\int_{s}^{t}(\frac{1}{\e}\cL_0\chi_{F})(X_{r}^{\e},\sL^{\mP}_{X_{r}^{\e}},Z_{r}^{\e})\dif r
+\int_{s}^{t}\frac{\e\partial \chi_{F}}{\partial x_{i}}(X_{r}^{\e},\sL^{\mP}_{X_{r}^{\e}},Z_{r}^{\e})\s_{1}^{ij}(X_{r}^{\e},\sL^{\mP}_{X_{r}^{\e}},Z_{r}^{\e})\dif B_{r}^{j}\\
&&+\int_{s}^{t}\frac{\e\partial \chi_{F}}{\partial z_{i}}(X_{r}^{\e},\sL^{\mP}_{X_{r}^{\e}},Z_{r}^{\e})\s_{2}^{ik}(X_{r}^{\e},\sL^{\mP}_{X_{r}^{\e}},Z_{r}^{\e})\dif W_{r}^{k}.
\de
So, multiplying $\Gamma_{s}$ and taking the expectation on both sides of the above equality, we know that
\ce
&&\mE\Bigg[\Gamma_{s}(\bar{X}_{.})\(\e\chi_{F}(X_{t}^{\e},\sL^{\mP}_{X_{t}^{\e}},Z_{t}^{\e})-\e\chi_{F}(X_{s}^{\e},\sL^{\mP}_{X_{s}^{\e}},Z_{s}^{\e})
-\e\int_{s}^{t}(\cL\chi_{F})(X_{r}^{\e},\sL^{\mP}_{X_{r}^{\e}},Z_{r}^{\e})\dif r\)\Bigg]\\
&=&\mE\Bigg[\Gamma_{s}(\bar{X}_{.})\int_{s}^{t}(\cL_0\chi_{F})(X_{r}^{\e},\sL^{\mP}_{X_{r}^{\e}},Z_{r}^{\e})\dif r\Bigg],
\de
which together with (\ref{limer}) yields that
\ce
&&\e\mE\Bigg[\Gamma_{s}(\bar{X}_{.})\(\chi_{F}(X_{t}^{\e},\sL^{\mP}_{X_{t}^{\e}},Z_{t}^{\e})-\chi_{F}(X_{s}^{\e},\sL^{\mP}_{X_{s}^{\e}},Z_{s}^{\e})-\int_{s}^{t}(\cL\chi_{F})(X_{r}^{\e},\sL^{\mP}_{X_{r}^{\e}},Z_{r}^{\e})\dif r\)\Bigg]\\
&=&-\mE\[\Gamma_{s}(\bar{X}_{.})\(F(X_{t}^{\e},\sL^{\mP}_{X_{t}^{\e}})-F(X_{s}^{\e},\sL^{\mP}_{X_{s}^{\e}})-\int_{s}^{t}(\bar{\cL}F)(X_{r}^{\e},\sL^{\mP}_{X_{r}^{\e}})\dif r\)\].
\de
Taking the limits on two sides of the above equality, by the boundedness of $\Gamma_{s}(\bar{X}_{.})$ and Proposition \ref{chib} we obtain that
\ce
\lim\limits_{\e\rightarrow0}\mE\[\Gamma_{s}(\bar{X}_{.})\(F(X_{t}^{\e},\sL^{\mP}_{X_{t}^{\e}})-F(X_{s}^{\e},\sL^{\mP}_{X_{s}^{\e}})
-\int_{s}^{t}(\bar{\cL}F)(X_{r}^{\e},\sL^{\mP}_{X_{r}^{\e}})\dif r\)\]=0,
\de
which implies that the weak limit of $\{X_{t}^{\e}, t\in[0,T]\}$ is a solution of the martingale problem associated with $(\bar{\cL},\d_{x_0})$. Since the solution of the martingale problem associated with $(\bar{\cL},\d_{x_0})$ is unique, $\{X_{t}^{\e}, t\in[0,T]\}$ converges weakly to $\{\bar{X}_{t}, t\in[0,T]\}$, which completes the proof.
\end{proof}

Now, it is the position to prove Theorem \ref{weakconv}.

{\bf Proof of Theorem \ref{weakconv}.}

Proposition \ref{xeweakcomp} and \ref{xeconx} yield Theorem \ref{weakconv}.

\section{Proof of Theorem \ref{convfilt}}\label{proofilt}

In this section, we prove Theorem \ref{convfilt}. The proof consists of three part. In the first part (Subsection \ref{sopres}), we present some preliminary estimates for $\Lambda^{\e}, \bar{\Lambda}, \rho_{t}^{\e}(1), \bar{\rho}_{t}(1)$. Then, we prove that $\{\xi_t^{\e}:=\rho_t^{\e,x,\mu}-\bar{\rho}_t,  t\in[0,T]\}$ is relatively weakly compact in the second part (Subsection \ref{rwcxi}). In the third part (Subsection \ref{weconze}), we establish that $\{\xi^{\e_k}\}$ converges weakly to $0$.

\subsection{Some preliminary estimates for $\Lambda^{\e}, \bar{\Lambda}, \rho_{t}^{\e}(1), \bar{\rho}_{t}(1)$}\label{sopres}

\bl\label{lambesti}
Suppose that $(\mathbf{H}^1_{b_{1}, \s_{1}})$, $(\mathbf{H}^{1}_{b_{2}, \s_{2}})$ and $(\mathbf{H}_{h})$ hold. Then we have that for any $q\geq1$ and $s,t\in[0,T]$,
\ce
\mE^{\mP^{\e}}|\Lambda_{t}^{\e}|^{2q}\leq C, \qquad \mE^{\mP^{\e}}|\Lambda_{s}^{\e}-\Lambda_{t}^{\e}|^{2q}\leq C|s-t|^q,
\de
where the constant $C$ is independent of $\e$.
\el
\begin{proof}
First of all, note that
\ce
\Lambda_{t}^{\e}=\exp\left\{\int_{0}^{t}h^{i}(X_{s}^{\e},\sL^{\mP}_{X_{s}^{\e}},Z_{s}^{\e})\dif Y_{s}^{\e,i} -\frac{1}{2}\int_{0}^{t}|h(X_{s}^{\e},\sL^{\mP}_{X_{s}^{\e}},Z_{s}^{\e})|^{2}\dif s\right\}.
\de
Thus, by the It\^{o} formula, it holds that
\be
\Lambda_{t}^{\e}=1+\int_{0}^{t}\Lambda_{s}^{\e}h^{i}(X_{s}^{\e},\sL^{\mP}_{X_{s}^{\e}},Z_{s}^{\e})\dif Y_{s}^{\e,i}.
\label{lamb}
\ee
Moreover, based on the BDG inequality  and $(\mathbf{H}_{h})$, one could obtain that
\ce
\mE^{\mP^{\e}}|\Lambda_{t}^{\e}|^{2q}
&\leq&2^{2q-1}+2^{2q-1}l^{2q-1}\sum_{i=1}^l\mE^{\mP^{\e}}\(\int_{0}^{t}\Lambda_{s}^{\e}h^{i}(X_{s}^{\e},\sL^{\mP}_{X_{s}^{\e}},Z_{s}^{\e})\dif Y_{s}^{\e,i}\)^{2q}\\
&\leq&2^{2q-1}+C\sum_{i=1}^l\int_{0}^{t}\mE^{\mP^{\e}}|\Lambda_{s}^{\e}h^i(X_{s}^{\e},\sL^{\mP}_{X_{s}^{\e}},Z_{s}^{\e})|^{2q}\dif s\\
&\leq&2^{2q-1}+C\int_{0}^{t}\mE^{\mP^{\e}}|\Lambda_{s}^{\e}|^{2q}\dif s,
\de
which together with the Gronwall inequality yields that
\be
\mE^{\mP^{\e}}(\Lambda_{t}^{\e})^{2q}\leq C.
\label{late}
\ee

Finally, we investigate that for $0\leq t<s\leq T$,
\ce
\Lambda_{s}^{\e}-\Lambda_{t}^{\e}=\int_{t}^{s}\Lambda_{r}^{\e}h^{i}(X_{r}^{\e},\sL^{\mP}_{X_{r}^{\e}},Z_{r}^{\e})\dif Y_{r}^{\e,i},
\de
and furthermore
\ce
\mE^{\mP^{\e}}|\Lambda_{s}^{\e}-\Lambda_{t}^{\e}|^{2q}
&\leq&l^{2q-1}\sum_{i=1}^l\mE^{\mP^{\e}}\(\int_{t}^{s}\Lambda_{r}^{\e}h^{i}(X_{r}^{\e},\sL^{\mP}_{X_{r}^{\e}},Z_{r}^{\e})\dif Y_{r}^{\e,i}\)^{2q}\\
&\leq&C(s-t)^{q-1}\sum_{i=1}^l\int_{t}^{s}\mE^{\mP^{\e}}|\Lambda_{r}^{\e}h^i(X_{r}^{\e},\sL^{\mP}_{X_{r}^{\e}},Z_{r}^{\e})|^{2q}\dif r\\
&\leq&C(s-t)^{q-1}\sum_{i=1}^l\int_{t}^{s}\mE^{\mP^{\e}}|\Lambda_{r}^{\e}|^{2q}\dif r\\
&\leq&C(s-t)^{q},
\de
where the last inequality is based on (\ref{late}). 
\end{proof}

\bl\label{-2es}
Under $(\mathbf{H}^1_{b_{1}, \s_{1}})$, $(\mathbf{H}^{1}_{b_{2}, \s_{2}})$ and $(\mathbf{H}_{h})$, it holds that
\ce
\sup\limits_{t\in[0,T]}\mE|\Lambda_{t}^{\e}|^{-2}\leq C.
\de
\el
\begin{proof}
Note that
\ce
\mE(\Lambda_{t}^{\e})^{-2}
&=&\mE\left[\left(\exp\left\{-\int_{0}^{t}h^{i}(X_{s}^{\e},\sL^{\mP}_{X_{s}^{\e}},Z_{s}^{\e})\dif V_{s}^{i} -\frac{1}{2}\int_{0}^{t}|h(X_{s}^{\e},\sL^{\mP}_{X_{s}^{\e}},Z_{s}^{\e})|^{2}\dif s\right\}\right)^{2}\right]\no\\
&=&\mE\Bigg[\exp\left\{-\int_{0}^{t}2 h^{i}(X_{s}^{\e},\sL^{\mP}_{X_{s}^{\e}},Z_{s}^{\e})\dif V_{s}^{i} -\frac{1}{2}\int_{0}^{t}|2 h(X_{s}^{\e},\sL^{\mP}_{X_{s}^{\e}},Z_{s}^{\e})|^{2}\dif s\right\}\no\\
&&\quad\cdot \exp\left\{\frac{1}{2}\int_{0}^{t}|2h(X_{s}^{\e},\sL^{\mP}_{X_{s}^{\e}},Z_{s}^{\e})|^{2}\dif s-\int_{0}^{t}|h(X_{s}^{\e},\sL^{\mP}_{X_{s}^{\e}},Z_{s}^{\e})|^{2}\dif s\right\}\Bigg]\no\\
&\leq& \exp\{CT\}\mE\Bigg[\exp\left\{-\int_{0}^{t}2 h^{i}(X_{s}^{\e},\sL^{\mP}_{X_{s}^{\e}},Z_{s}^{\e})\dif V_{s}^{i} -\frac{1}{2}\int_{0}^{t}|2 h(X_{s}^{\e},\sL^{\mP}_{X_{s}^{\e}},Z_{s}^{\e})|^{2}\dif s\right\}\Bigg]\no\\
&=&\exp\{CT\},
\de
where the last equality is based on the fact that $\exp\big\{-\int_{0}^{t}2 h^{i}(X_{s}^{\e},\sL^{\mP}_{X_{s}^{\e}},Z_{s}^{\e})\dif V_{s}^{i}\\ -\frac{1}{2}\int_{0}^{t}|2 h(X_{s}^{\e},\sL^{\mP}_{X_{s}^{\e}},Z_{s}^{\e})|^{2}\dif s\big\}$ is an exponential martingale. The proof is complete.
\end{proof}

\bl\label{lambbaresti}
Under $(\mathbf{H}^1_{b_{1}, \s_{1}})$, $(\mathbf{H}^2_{\s_{1}})$, $(\mathbf{H}^{1}_{b_{2}, \s_{2}})$-$(\mathbf{H}^{2}_{b_{2}, \s_{2}})$, $(\mathbf{H}_{h})$, there exists a constant $C>0$ such that for any $q\geq1$,
\ce
\sup\limits_{t\in[0,T]}\mE^{\mP^{\e}}|\bar{\Lambda}_{t}|^{2q}\leq C.
\de
\el

Since the proof of the above lemma is similar to that for $\Lambda^\e$ in Lemma \ref{lambesti}, we omit it.

Next, in order to estimate $\rho_{t}^{\e}(1), \bar{\rho}_{t}(1)$, we establish the Zakai equations about $\rho^\e, \bar{\rho}$.

\bl  (The Zakai equation)  \label{zaki}\\
(i) For $\Psi\in C^{2,(1,1),2}_{b}(\mR^{n}\times\cP_{2}(\mR^{n})\times\mR^{m})$, the Zakai equation of Eq.(\ref{slfa}) is given by
\be
\rho_{t}^{\e}(\Psi)=\rho_{0}^{\e}(\Psi)+\int_{0}^{t}\rho_{s}^{\e}\((\cL+\frac{1}{\e}\cL_0)\Psi\)\dif s
+\int_{0}^{t}\rho_{s}^{\e}(\Psi h^{i})\dif Y_{s}^{\e,i},\quad \rho_{0}^{\e}(\Psi)=\Psi(x_0,\d_{x_0},z_0),
\label{rhvi}
\ee
(ii) For $F\in C^{2,(1,1)}_{b}(\mR^{n}\times\cP_{2}(\mR^{n}))$, the Zakai equation of Eq.(\ref{Eq3}) is given by
\be
\bar{\rho}_{t}(F)=\bar{\rho}_{0}(F)+\int_{0}^{t}\bar{\rho}_{s}(\bar{\cL}F)\dif s
+\int_{0}^{t}\bar{\rho}_{s}(F \bar{h}^{i})\dif Y_{s}^{\e,i},\quad \bar{\rho}_{0}(F)=F(x_0,\d_{x_0}).
\label{barh}
\ee
\el
\begin{proof}
Applying the It\^{o} formula to $\Psi(X_{t}^{\e},\sL^{\mP}_{X_{t}^{\e}},Z_{t}^{\e})$, one could obtain that
\ce
\Psi(X_{t}^{\e},\sL^{\mP}_{X_{t}^{\e}},Z_{t}^{\e})&=&\Psi(x_0,\d_{x_0},z_0)
+\int_{0}^{t}\((\cL+\frac{1}{\e}\cL_0)\Psi\)(X_{s}^{\e},\sL^{\mP}_{X_{s}^{\e}},Z_{s}^{\e})\dif s\\
&&+\int_{0}^{t}\frac{\partial \Psi}{\partial x_{i}}(X_{s}^{\e},\sL^{\mP}_{X_{s}^{\e}},Z_{s}^{\e})\s_{1}^{ij}(X_{s}^{\e},\sL^{\mP}_{X_{s}^{\e}},Z_{s}^{\e})\dif B_{s}^{j}\\
&&+\int_{0}^{t}\frac{\partial \Psi}{\partial z_{i}}(X_{s}^{\e},\sL^{\mP}_{X_{s}^{\e}},Z_{s}^{\e})\s_{2}^{ik}(X_{s}^{\e},\sL^{\mP}_{X_{s}^{\e}},Z_{s}^{\e})\dif W_{s}^{k}.
\de
So, by combining the above deduction with (\ref{lamb}), the It\^{o} formula implies that
\ce
\Psi(X_{t}^{\e},\sL^{\mP}_{X_{t}^{\e}},Z_{t}^{\e})\Lambda_{t}^{\e}
&=&\Psi(x_0,\d_{x_0},z_0)+\int_{0}^{t}\Lambda_{s}^{\e}\((\cL+\frac{1}{\e}\cL_0)\Psi\)(X_{s}^{\e},\sL^{\mP}_{X_{s}^{\e}},Z_{s}^{\e})\dif s\\
&&+\int_{0}^{t}\Psi(X_{s}^{\e},\sL^{\mP}_{X_{s}^{\e}},Z_{s}^{\e})\Lambda_{s}^{\e}h^{i}(X_{s}^{\e},\sL^{\mP}_{X_{s}^{\e}})\dif Y_{s}^{\e,i}\\
&&+\int_{0}^{t}\Lambda_{s}^{\e}\frac{\partial \Psi}{\partial x_{i}}(X_{s}^{\e},\sL^{\mP}_{X_{s}^{\e}},Z_{s}^{\e})\s_{1}^{ij}(X_{s}^{\e},\sL^{\mP}_{X_{s}^{\e}},Z_{s}^{\e})\dif B_{s}^{j}\\
&&+\int_{0}^{t}\Lambda_{s}^{\e}\frac{\partial \Psi}{\partial z_{i}}(X_{s}^{\e},\sL^{\mP}_{X_{s}^{\e}},Z_{s}^{\e})\s_{2}^{ik}(X_{s}^{\e},\sL^{\mP}_{X_{s}^{\e}},Z_{s}^{\e})\dif W_{s}^{k}.
\de
Taking the conditional expectation about  $\mathscr{F}_{t}^{Y^{\varepsilon}}$ on  both sides of the above equality under the measure $\mP^{\e}$, one can get that
\ce
&&\mE^{\mP^{\e}}[\Psi(X_{t}^{\e},\sL^{\mP}_{X_{t}^{\e}},Z_{t}^{\e})\Lambda_{t}^{\e}|\mathscr{F}_{t}^{Y^{\e}}]\\
&=&\Psi(x_0,\d_{x_0},z_0)
+\int_{0}^{t}\mE^{\mP^{\e}}\[\Lambda_{s}^{\e}\((\cL+\frac{1}{\e}\cL_0)\Psi\)(X_{s}^{\e},\sL^{\mP}_{X_{s}^{\e}},Z_{s}^{\e})|\mathscr{F}_{s}^{Y^{\e}}\]\dif s\\
&&+\int_{0}^{t}\mE^{\mP^{\e}}[\Psi(X_{s}^{\e},\sL^{\mP}_{X_{s}^{\e}},Z_{s}^{\e})\Lambda_{s}^{\e}h^{i}(X_{s}^{\e},\sL^{\mP}_{X_{s}^{\e}})|\mathscr{F}_{s}^{Y^{\e}}]\dif Y_{s}^{\e,i}\\
&=&\mE^{\mP^{\e}}[\Psi(x_0,\d_{x_0},z_0)|\mathscr{F}_{0}^{Y^{\e}}]
+\int_{0}^{t}\mE^{\mP^{\e}}\[\Lambda_{s}^{\e}\((\cL+\frac{1}{\e}\cL_0)\Psi\)(X_{s}^{\e},\sL^{\mP}_{X_{s}^{\e}},Z_{s}^{\e})|\mathscr{F}_{s}^{Y^{\e}}\]\dif s\\
&&+\int_{0}^{t}\mE^{\mP^{\e}}[\Psi(X_{s}^{\e},\sL^{\mP}_{X_{s}^{\e}},Z_{s}^{\e})\Lambda_{s}^{\e}h^{i}(X_{s}^{\e},\sL^{\mP}_{X_{s}^{\e}})|\mathscr{F}_{s}^{Y^{\e}}]\dif Y_{s}^{\e,i},
\de
which together with the definition of $\rho_{t}^{\e}(\Psi)$ yields that
\ce
\rho_{t}^{\e}(\Psi)=\rho_{0}^{\e}(\Psi)+\int_{0}^{t}\rho_{s}^{\e}\((\cL+\frac{1}{\e}\cL_0)\Psi\)\dif s
+\int_{0}^{t}\rho_{s}^{\e}(\Psi h^{i})\dif Y_{s}^{\e,i}, \quad \rho_{0}^{\e}(\Psi)=\Psi(x_0,\d_{x_0},z_0).
\de

By the same deduction to that of $(i)$, we obtain $(ii)$. The proof is complete.
\end{proof}

\bl\label{rhbrh1}
Under $(\mathbf{H}^1_{b_{1}, \s_{1}})$, $(\mathbf{H}^2_{\s_{1}})$, $(\mathbf{H}^{1}_{b_{2}, \s_{2}})$-$(\mathbf{H}^{2}_{b_{2}, \s_{2}})$, $(\mathbf{H}_{h})$, there exists a constant $C>0$ such that for any $q\geq1$,
\ce
\mE\left(\sup\limits_{t\in[0,T]}|\rho_{t}^{\e}(1)|^{q}\right)\leq C, \quad\mE\left(\sup\limits_{t\in[0,T]}|\bar{\rho}_{t}(1)|^{q}\right)\leq C.
\de
\el
\begin{proof}
First of all, from the H\"{o}lder inequality, it follows that
\be
\mE\left(\sup\limits_{t\in[0,T]}|\rho_{t}^{\e}(1)|^{q}\right)=\mE^{\mP^{\e}}\left[\Lambda_{T}^{\e}\left(\sup\limits_{t\in[0,T]}|\rho_{t}^{\e}(1)|^{q}\right)\right]
=[\mE^{\mP^{\e}}(\Lambda_{T}^{\e})^{2}]^{\frac{1}{2}}\left[\mE^{\mP^{\e}}\left(\sup\limits_{t\in[0,T]}|\rho_{t}^{\e}(1)|^{2q}\right)\right]^{\frac{1}{2}}.
\label{rhop}
\ee
Then we estimate $|\rho_{t}^{\e}(1)|^{2q}$. Since $\rho_{t}^{\e}(1)$ satisfies Eq.(\ref{rhvi}) with $\Psi(x,\mu,z)=1$, it holds that
\ce
\mE^{\mP^{\e}}\left(\sup\limits_{t\in[0,T]}|\rho_{t}^{\e}(1)|^{2q}\right)
&\leq&2^{2q-1}+2^{2q-1}l^{2q-1}\sum_{i=1}^l\mE^{\mP^{\e}}\left(\sup\limits_{t\in[0,T]}\Big|\int_{0}^{t}\rho_{s}^{\e}(h^{i})\dif Y_{s}^{\e,i}\Big|^{2q}\right)\no\\
&\leq&2^{2q-1}+2^{2q-1}l^{2q-1}\Big[\frac{2q(2q-1)}{2}\Big]^{q}\sum_{i=1}^l\mE^{\mP^{\e}}\(\int_{0}^{T}|\rho_{s}^{\e}(h^{i})|^2\dif s\)^{q}\no\\
&\leq& 2^{2q-1}+2^{2q-1}l^{2q-1}\Big[\frac{2q(2q-1)}{2}\Big]^{q}T^{q-1}\sum_{i=1}^l\int_{0}^{T}\mE^{\mP^{\e}}|\rho_{s}^{\e}(h^{i})|^{2q}\dif s\no\\
&\leq&2^{2q-1}+C\sum_{i=1}^l\int_{0}^{T}\mE^{\mP^{\e}}\Big[\Big|\mE^{\mP^{\e}}[h^i(X_{s}^{\e},\sL^{\mP}_{X_{s}^{\e}},Z_{s}^{\e})\Lambda_{s}^{\e}|\mathscr{F}_{s}^{Y^{\e}}]\Big|^{2q}\Big]\dif s\no\\
&\leq&2^{2q-1}+C\sum_{i=1}^l\int_{0}^{T}\mE^{\mP^{\e}}\Big|h^i(X_{s}^{\e},\sL^{\mP}_{X_{s}^{\e}},Z_{s}^{\e})\Lambda_{s}^{\e}\Big|^{2q}\dif s\no\\
&\leq& 2^{2q-1}+C\int_{0}^{T}\mE^{\mP^{\e}}|\Lambda_{s}^{\e}|^{2q}\dif s\\
&\leq& C,
\de
where the last step is based on (\ref{late}). So, inserting the above inequality in (\ref{rhop}), by (\ref{late}) we obtain that
$$
\mE\left(\sup\limits_{t\in[0,T]}|\rho_{t}^{\e}(1)|^{q}\right)\leq C.
$$

Finally, the same deduction to that for the above inequality implies the second estimate. The proof is complete.
\end{proof}

\bl\label{rhbi}
Under $(\mathbf{H}^1_{b_{1}, \s_{1}})$, $(\mathbf{H}^2_{\s_{1}})$, $(\mathbf{H}^{1}_{b_{2}, \s_{2}})$-$(\mathbf{H}^{2}_{b_{2}, \s_{2}})$, $(\mathbf{H}_{h})$, there exists a constant $C>0$ such that for $t\in[0,T]$
\ce
(\bar{\rho}_{t}(1))^{-1}\leq C, \quad a.s..
\de
\el

Since the proof of the above lemma is similar to that of \cite[Lemma 5.1]{Qiao1}, we omit it.

\subsection{Relatively weak compactness for $\{\xi_t^{\e}:=\rho_t^{\e,x,\mu}-\bar{\rho}_t,  t\in[0,T]\}$}\label{rwcxi} 

In order to investigate $\xi^{\e}$, we prepare some following stronger moment estimates.

\bl \label{xezebp}
 Under $(\mathbf{H}^{1}_{b_{1}, \s_{1}})$ $(\mathbf{H}^{1}_{b_{2}, \s_{2}})$-$(\mathbf{H}^{2'}_{b_{2}, \s_{2}})$, there exists a constant $C>0$ such that 
\ce
\sup\limits_{t\in[0,T]}\mE|X_{t}^{\e}|^{2p}\leq C(1+|x_0|^{2p}+|z_0|^{2p}), \quad \sup\limits_{t\in[0,T]}\mE|Z_{t}^{\e}|^{2p}\leq C(1+|x_0|^{2p}+|z_0|^{2p}),
\de
where $p$ is the same to that in $(\mathbf{H}^{2'}_{b_{2}, \s_{2}})$.
\el

 We omit the proof of the above lemma, since it is similar to that for Lemma \ref{xezeb}.

\bp\label{corbr}
Under $(\mathbf{H}^1_{b_{1}, \s_{1}})$ $(\mathbf{H}^2_{\s_{1}})$ $(\mathbf{H}^{1}_{b_{2}, \s_{2}})$-$(\mathbf{H}^{2'}_{b_{2}, \s_{2}})$ $(p\geq 8)$,  and $(\mathbf{H}_{h})$, $\{\xi_{t}^{\e}, t\in[0,T]\}$  is relatively weakly compact in $ C([0,T],\cM(\mR^{n}\times\sP_{2}(\mR^{n})))$, where $\cM(\mR^{n}\times\sP_{2}(\mR^{n}))$ stands for the collection of all the finite measures on $\mR^{n}\times\sP_{2}(\mR^{n})$.
\ep
\begin{proof}
First of all, from Lemma \ref{rhbrh1}, it follows that for $F\in C^{2,(1,1)}_{b}(\mR^{n}\times\cP_{2}(\mR^{n}))$,
\ce
\mE|\xi_{t}^{\e}(F)|&=&\mE|\rho_{t}^{\e,x,\mu}(F)-\bar{\rho}_{t}(F)|=\mE|\rho_{t}^{\e}(F)-\bar{\rho}_{t}(F)|
\leq\mE|\rho_{t}^{\e}(F)|+\mE|\bar{\rho}_{t}(F)|\\
&\leq&\|F\|_{C^{2,(1,1)}_{b}(\mR^{n}\times\cP_{2}(\mR^{n}))}\mE|\rho_{t}^{\e}(1)|+\|F\|_{C^{2,(1,1)}_{b}(\mR^{n}\times\cP_{2}(\mR^{n}))}\mE|\bar{\rho}_{t}(1)|\\
&\leq&C\|F\|_{C^{2,(1,1)}_{b}(\mR^{n}\times\cP_{2}(\mR^{n}))},
\de
which yields that
\be
\sup_{\e}\sup_{t\in[0,T]}\mE|\xi_{t}^{\e}(F)|<\infty.
\label{surh}
\ee

Besides, by Lemma \ref{zaki}, it holds that
\ce
&&\rho_{t}^{\e,x,\mu}(F)=F(x_0,\d_{x_0})+\int_{0}^{t}\rho_{s}^{\e,x,\mu}\left((\cL F)(\cdot,\cdot,Z_s^\e)\right)
\dif s+\int_{0}^{t}\rho_{s}^{\e,x,\mu}\left(F h^{i}(\cdot,\cdot,Z_s^\e)\right)\dif Y_{s}^{\e,i},\\
&&\bar{\rho}_{t}(F)=F(x_0,\d_{x_0})+\int_{0}^{t}\bar{\rho}_{s}(\bar{\cL}F)\dif s
+\int_{0}^{t}\bar{\rho}_{s}(F \bar{h}^{i})\dif Y_{s}^{\e,i}.
\de
Based on the above equation and the isometric formula, we can obtain that for any $\delta>0$ and any $\{\mathscr{F}_{t}\}$-stopping time $0\leq\tau<\tau+\d\leq T$,
\ce
\mE|\xi_{\tau+\d}^{\e}(F)-\xi_{\tau}^{\e}(F)|^{2}&=&\mE^{\mP^{\e}}\left[|\xi_{\tau+\d}^{\e}(F)-\xi_{\tau}^{\e}(F)|^{2}\Lambda_{T}^{\e}\right]\\
&\leq&(\mE^{\mP^{\e}}|\xi_{\tau+\d}^{\e}(F)-\xi_{\tau}^{\e}(F)|^{4})^{1/2}(\mE^{\mP^{\e}}(\Lambda_{T}^{\e})^2)^{1/2},
\de
and
\ce
\mE^{\mP^{\e}}|\xi_{\tau+\d}^{\e}(F)-\xi_{\tau}^{\e}(F)|^{4}&\leq&\mE^{\mP^{\e}}|(\rho_{\tau+\d}^{\e,x,\mu}(F)-\bar{\rho}_{\tau+\d}(F))-(\rho_{\tau}^{\e,x,\mu}(F)-\bar{\rho}_{\tau}(F))|^{4}\\
&\leq&2^3\mE^{\mP^{\e}}|\rho_{\tau+\d}^{\e,x,\mu}(F)-\rho_{\tau}^{\e,x,\mu}(F)|^{4}+2^3\mE^{\mP^{\e}}|\bar{\rho}_{\tau+\d}(F)-\bar{\rho}_{\tau}(F)|^{4}\\
&\leq&4^3\mE^{\mP^{\e}}\Big|\int_{\tau}^{\tau+\d}\rho_{s}^{\e,x,\mu}\left((\cL F)(\cdot,\cdot,Z_s^\e)\right)\dif s\Big|^{4}\\
&&+4^3\mE^{\mP^{\e}}\Big|\int_{\tau}^{\tau+\d}\rho_{s}^{\e,x,\mu}\left(F h^{i}(\cdot,\cdot,Z_s^\e)\right)\dif Y_{s}^{\e,i}\Big|^{4}\\
&&+4^3\mE^{\mP^{\e}}\Big|\int_{\tau}^{\tau+\d}\bar{\rho}_{s}(\bar{\cL}F)\dif s\Big|^{4}
+4^3\mE^{\mP^{\e}}\Big|\int_{\tau}^{\tau+\d}\bar{\rho}_{s}(F \bar{h}^{i})\dif Y_{s}^{\e,i}\Big|^{4}\\
&\leq&4^3\d^3\mE^{\mP^{\e}}\int_{\tau}^{\tau+\d}\left|\rho_{s}^{\e,x,\mu}\left((\cL F)(\cdot,\cdot,Z_s^\e)\right)\right|^{4}\dif s\\
&&+4^3\d\mE^{\mP^{\e}}\int_{\tau}^{\tau+\d}\left|\rho_{s}^{\e,x,\mu}\left(F h^{i}(\cdot,\cdot,Z_s^\e)\right)\right|^{4}\dif s\\
&&+4^3\d^3\mE^{\mP^{\e}}\int_{\tau}^{\tau+\d}|\bar{\rho}_{s}(\bar{\cL}F)|^{4}\dif s
+4^3\d\mE^{\mP^{\e}}\int_{\tau}^{\tau+\d}|\bar{\rho}_{s}(F \bar{h}^{i})|^{4}\dif s\\
&=:&I_{1}+I_{2}+I_{3}+I_{4}.
\de

For $I_{1}+I_3$, it holds that
\ce
I_{1}+I_3&\leq&4^3\d^3\mE^{\mP^{\e}}\int_{\tau}^{\tau+\d}\mE^{\mP^{\e}}\[|(\cL F)(X_s^{\e},\sL^{\mP}_{X_s^{\e}},Z_s^\e)|^{4}|\Lambda^\e_{s}|^4|\mathscr{F}_{s}^{Y^{\e}}\]\dif s\\
&&+4^3\d^3\mE^{\mP^{\e}}\int_{\tau}^{\tau+\d}\mE^{\mP^{\e}}\[|(\bar{\cL} F)(\bar{X}_{s},\sL^{\mP}_{\bar{X}_{s}})|^4|\bar{\Lambda}_{s}|^{4}|\mathscr{F}_{s}^{Y^{\e}}\]\dif s\\
&\leq&4^3\d^3 C\mE^{\mP^{\e}}\int_{\tau}^{\tau+\d}\mE^{\mP^{\e}}\[(1+|X_s^{\e}|^8+\mE|X_s^{\e}|^8+|Z_s^\e|^{8})|\Lambda^\e_{s}|^4|\mathscr{F}_{s}^{Y^{\e}}\]\dif s\\
&&+4^3\d^3 C\mE^{\mP^{\e}}\int_{\tau}^{\tau+\d}\mE^{\mP^{\e}}\[(1+|\bar{X}_{s}|^8+\mE|\bar{X}_{s}|^8)|\bar{\Lambda}_{s}|^{4}|\mathscr{F}_{s}^{Y^{\e}}\]\dif s\\
&\leq&4^3\d^3 C\mE^{\mP^{\e}}\int_{0}^{\d}\mE^{\mP^{\e}}\[(1+|X_{\tau+u}^{\e}|^8+C+|Z_{\tau+u}^\e|^{8})|\Lambda^\e_{{\tau+u}}|^4|\mathscr{F}_{{\tau+u}}^{Y^{\e}}\]\dif u\\
&&+4^3\d^3 C\mE^{\mP^{\e}}\int_{0}^{\d}\mE^{\mP^{\e}}\[(1+|\bar{X}_{{\tau+u}}|^8+C)|\bar{\Lambda}_{{\tau+u}}|^{4}|\mathscr{F}_{{\tau+u}}^{Y^{\e}}\]\dif u\\
&\leq&4^3\d^3 C\int_{0}^{\d}\left(1+\mE^{\mP^{\e}}|X_{\tau+u}^{\e}|^{16}+C+\mE^{\mP^{\e}}|Z_{\tau+u}^\e|^{16}\right)^{1/2}\left(\mE^{\mP^{\e}}|\Lambda_{\tau+u}^{\e}|^8\right)^{1/2}\dif u\\
&&+4^3\d^3 C\int_{0}^{\d}\left(1+\mE^{\mP^{\e}}|\bar{X}_{{\tau+u}}|^{16}+C\right)^{1/2}\left(\mE^{\mP^{\e}}|\bar{\Lambda}_{t}|^{8}\right)^{1/2}\dif u\\
&\leq& C\d^{4},
\de
where the last inequality is based on Lemma \ref{xezebp}, \ref{lambesti}, \ref{barsig} and \ref{lambbaresti}. Moreover, by the boundedness of $F, h$, we have that $I_{2}+I_{4}\leq C\d^2$. Thus, these estimates imply that
$$
\mE|\xi_{\tau+\d}^{\e}(F)-\xi_{\tau}^{\e}(F)|^{2}\leq C(\d^{2}+\d),
$$
and
\be
\lim_{\d\rightarrow0}\limsup_{\e\downarrow 0}\sup\limits_{\tau\leq T}\mE|\xi_{\tau+\d}^{\e}(F)-\xi_{\tau}^{\e}(F)|^{2}=0.
\label{lirh}
\ee

Finally, combining (\ref{surh}) (\ref{lirh}) with \cite[Theorem 2.7]{TGK}, we obtain that $\{\xi_{t}^{\e}(F), t\in[0,T]\}$ is relatively weakly compact in $C([0,T],\mR)$, which together with \cite[Theorem 6.2]{kus} yields that $\{\xi_{t}^{\e},t\in[0,T]\}$ is relatively weakly compact.
\end{proof}

\subsection{Weak convergence for $\{\xi_t^{\e},  t\in[0,T]\}$ to $0$}\label{weconze}

Here, in order to prove that $\{\xi^{\e_k}\}$ converges weakly to $0$, we need the following estimate.

\bl \label{xediffp}
 Suppose that $(\mathbf{H}^{1}_{b_{1}, \s_{1}})$ $(\mathbf{H}^{1}_{b_{2}, \s_{2}})$-$(\mathbf{H}^{2'}_{b_{2}, \s_{2}})$ hold. Then it holds that for any $\delta>0$ and any $t\in[0,T]$,
 \ce
 \mE|X_{t+\d}^{\e}-X_{t}^{\e}|^{2p}\leq C(1+|x_0|^{2p}+|z_0|^{2p})(\d^{2p}+\d^p),
 \de
 where $p$ is the same to that in $(\mathbf{H}^{2'}_{b_{2}, \s_{2}})$.
 \el
 
 We omit the proof of the above lemma, since it is similar to that for Lemma \ref{xediff}.

\bp\label{corbr2}
Under $(\mathbf{H}^1_{b_{1}, \s_{1}})$ $(\mathbf{H}^2_{\s_{1}})$ $(\mathbf{H}^3_{b_{1}, \s_{1}})$ $(\mathbf{H}^{1}_{b_{2}, \s_{2}})$ $(\mathbf{H}^{2'}_{b_{2}, \s_{2}}) (p\geq 12)$ $(\mathbf{H}^{3}_{b_{2}, \s_{2}})$ and $(\mathbf{H}_{h})$, there exists a subsequence $\{\xi^{\e_k}\}$ which converges weakly to $0$ in $ C([0,T],\cM(\mR^{n}\times\sP_{2}(\mR^{n})))$.
\ep
\begin{proof}
{\bf Step 1.} We prove that for any $t\in[0,T]$ and $F\in \mC^{4,(2,2)}_{b}(\mR^{n}\times\cP_{2}(\mR^{n}))$
\ce
\lim\limits_{\e\rightarrow 0}\mE\left|\xi_{t}^{\e}(F)-\int_{0}^{t}\xi_{s}^{\e}(\bar{\cL}F)\dif s-\int_{0}^{t}\xi_{s}^{\e}(F \bar{h}^{i})\dif Y_{s}^{\e,i}\right|^{2}=0.
\de

First of all, we define a perturbed test function $F^{\e}(x,\mu,z)$ as follows:
$$
F^{\e}(x,\mu,z)=F(x,\mu)+\e\chi_F(x,\mu,z),
$$
where $\chi_F(x,\mu,z)$ is the unique solution of the Poisson equation (\ref{pois}). Note that
$$
\rho_{t}^{\e}(F^{\e})=\rho_{t}^{\e}(F)+\rho_{t}^{\e}(\e\chi_F)=\rho_{t}^{\e,x,\mu}(F)+\e\rho_{t}^{\e}(\chi_F).
$$
Thus, based on Lemma \ref{zaki}, we get that
\ce
\xi_{t}^{\e}(F)&=&\rho_{t}^{\e,x,\mu}(F)-\bar{\rho}_{t}(F)=\rho_{t}^{\e}(F^{\e})-\e\rho_{t}^{\e}(\chi_F)-\bar{\rho}_{t}(F)\\
&\overset{(\ref{rhvi})(\ref{barh})}{=}&\rho_{0}^{\e}(F^{\e})+\int_{0}^{t}\rho_{s}^{\e}\((\cL+\frac{1}{\e}\cL_0)F^{\e}\)\dif s
+\int_{0}^{t}\rho_{s}^{\e}(F^{\e} h^{i})\dif Y_{s}^{\e,i}-\e\rho_{t}^{\e}(\chi_F)\\
&&-\(\bar{\rho}_{0}(F)+\int_{0}^{t}\bar{\rho}_{s}(\bar{\cL}F)\dif s
+\int_{0}^{t}\bar{\rho}_{s}(F \bar{h}^{i})\dif Y_{s}^{\e,i}\)\\
&=&\rho_{0}^{\e,x,\mu}(F)+\e\rho_{0}^{\e}(\chi_F)
+\int_{0}^{t}\rho_{s}^{\e}(\cL F+\e\cL\chi_F+\cL_0\chi_F)\dif s\\
&&+\int_{0}^{t}\rho_{s}^{\e}(F h^{i})\dif Y_{s}^{\e,i}+\e\int_{0}^{t}\rho_{s}^{\e}(\chi_F h^{i})\dif Y_{s}^{\e,i}-\e\rho_{t}^{\e}(\chi_F)
-\bar{\rho}_{0}(F)-\int_{0}^{t}\bar{\rho}_{s}(\bar{\cL}F)\dif s\\
&&-\int_{0}^{t}\bar{\rho}_{s}(F \bar{h}^{i})\dif Y_{s}^{\e,i}\\
&\overset{(\ref{pois})}{=}&\e\rho_{0}^{\e}(\chi_F)+\int_{0}^{t}\rho_{s}^{\e}(\bar{\cL}F)\dif s+\e\int_{0}^{t}\rho_{s}^{\e}(\cL\chi_F)\dif s\\
&&+\int_{0}^{t}\rho_{s}^{\e}(F h^{i})\dif Y_{s}^{\e,i}
+\e\int_{0}^{t}\rho_{s}^{\e}(\chi_F h^{i})\dif Y_{s}^{\e,i}
-\e\rho_{t}^{\e}(\chi_F)-\int_{0}^{t}\bar{\rho}_{s}(\bar{\cL}F)\dif s\\
&&-\int_{0}^{t}\bar{\rho}_{s}(F \bar{h}^{i})\dif Y_{s}^{\e,i}\\
&=&\e\rho_{0}^{\e}(\chi_F)
+\int_{0}^{t}\rho_{s}^{\e,x,\mu}(\bar{\cL}F)\dif s+\e\int_{0}^{t}\rho_{s}^{\e}(\cL\chi_F)\dif s+\int_{0}^{t}\rho_{s}^{\e}(F h^{i})\dif Y_{s}^{\e,i}\\
&&+\e\int_{0}^{t}\rho_{s}^{\e}(\chi_F h^{i})\dif Y_{s}^{\e,i}-\e\rho_{t}^{\e}(\chi_F)-\int_{0}^{t}\bar{\rho}_{s}(\bar{\cL}F)\dif s-\int_{0}^{t}\bar{\rho}_{s}(F \bar{h}^{i})\dif Y_{s}^{\e,i}\\
&&-\int_{0}^{t}\rho_{s}^{\e,x,\mu}(F \bar{h}^{i})\dif Y_{s}^{\e,i}+\int_{0}^{t}\rho_{s}^{\e,x,\mu}(F \bar{h}^{i})\dif Y_{s}^{\e,i}\\
&=&\e\rho_{0}^{\e}(\chi_F)+\int_{0}^{t}\xi_{s}^{\e}(\bar{\cL}F)\dif s+\e\int_{0}^{t}\rho_{s}^{\e}(\cL\chi_F)\dif s+\int_{0}^{t}\rho_{s}^{\e}(F h^{i}-F \bar{h}^{i})\dif Y_{s}^{\e,i}\\
&&+\e\int_{0}^{t}\rho_{s}^{\e}(\chi_F h^{i})\dif Y_{s}^{\e,i}-\e\rho_{t}^{\e}(\chi_F)+\int_{0}^{t}\xi_{s}^{\e}(F \bar{h}^{i})\dif Y_{s}^{\e,i},
\de
where the fact that $\rho_{0}^{\e,x,\mu}(F)=\bar{\rho}_{0}(F)=F(x_0, \d_{x_0})$ is used in the fourth equality. Moreover, it holds that
\be
&&\mE\left|\xi_{t}^{\e}(F)-\int_{0}^{t}\xi_{s}^{\e}(\bar{\cL}F)\dif s-\int_{0}^{t}\xi_{s}^{\e}(F \bar{h}^{i})\dif Y_{s}^{\e,i}\right|^{2}\no\\
&=&\mE\bigg|\e\rho_{0}^{\e}(\chi_F)-\e\rho_{t}^{\e}(\chi_F)+\e\int_{0}^{t}\rho_{s}^{\e}(\cL\chi_F)\dif s+\e\int_{0}^{t}\rho_{s}^{\e}(\chi_F h^{i})\dif Y_{s}^{\e,i}\no\\
&&\qquad +\int_{0}^{t}\rho_{s}^{\e}(F h^{i}-F \bar{h}^{i})\dif Y_{s}^{\e,i}\bigg|^{2}\no\\
&\leq&5\e^{2}\mE|\rho_{0}^{\e}(\chi_F)|^{2}+5\e^{2}\mE|\rho_{t}^{\e}(\chi_F)|^{2}+5\e^{2}\mE\left|\int_{0}^{t}\rho_{s}^{\e}(\cL\chi_F)\dif s\right|^{2}\no\\
&&+5\e^{2}\mE\left|\int_{0}^{t}\rho_{s}^{\e}(\chi_F h^{i})\dif Y_{s}^{\e,i}\right|^{2}+5\mE\left|\int_{0}^{t}\rho_{s}^{\e}(F h^{i}-F \bar{h}^{i})\dif Y_{s}^{\e,i}\right|^{2}\no\\
&=:&J_{1}+J_{2}+J_{3}+J_{4}+J_{5}.
\label{mxifc}
\ee

In the following we treat $J_1$. It is easy to see that
\be
J_1=5\e^{2}|\chi_F(x_0,\d_{x_0},z_0)|^{2}.
\label{j1}
\ee
For $J_{2}$, based on the H\"{o}lder inequality, the Jensen inequality, Lemma \ref{lambesti}, \ref{xezebp} and Proposition \ref{chib}, it holds that
\be
J_{2}&=& 5\e^{2}\mE^{\mP^{\e}}[|\rho_{t}^{\e}(\chi_F)|^{2}\Lambda_T^\e]\leq 5\e^{2}(\mE^{\mP^{\e}}|\rho_{t}^{\e}(\chi_F)|^{4})^{1/2}(\mE^{\mP^{\e}}|\Lambda_T^\e|^2)^{1/2}\no\\
&\leq& 5\e^{2}C(\mE^{\mP^{\e}}[\mE^{\mP^{\e}}[|\chi_F(X_{t}^{\e},\sL^{\mP}_{X^{\e}},Z_{t}^{\e})|^4|\Lambda_{t}^{\e}|^4|\mathscr{F}_{t}^{Y^{\e}}]])^{1/2}\no\\
&\leq& 5\e^{2}C(\mE^{\mP^{\e}}|\chi_F(X_{t}^{\e},\sL^{\mP}_{X^{\e}},Z_{t}^{\e})|^8)^{1/4}(\mE^{\mP^{\e}}|\Lambda_{t}^{\e}|^8)^{1/4}\no\\
&\leq& 5\e^{2}C(\mE^{\mP^{\e}}(1+\|\sL^{\mP}_{X^{\e}}\|^8+|Z_{t}^{\e}|^8))^{1/4}\no\\
&\leq& C\e^{2}.
\label{j2}
\ee
For $J_3$, Proposition \ref{chib} and (\ref{b1line}) imply that
$$
|(\cL\chi_F)(x,\mu,z)|\leq C(1+|x|+\|\mu\|+|z|)^3. 
$$
And based on the H\"{o}lder inequality and the same deduction to that for $J_{2}$, we know that
\be
J_{3}\leq5T\e^{2}\int_{0}^{T}\mE|\rho_{s}^{\e}(\cL\chi_F)|^{2}\dif s\leq C\e^{2}.
\label{j3}
\ee
For $J_{4}$, by the boundedness of $h$, the BDG inequality and the similar deduction to that for $J_{1}$, it holds that
\be
J_{4}&=&5\e^{2}\mE^{\mP^{\e}}\Bigg[\Lambda_{T}^{\e}\left|\int_{0}^{t}\rho_{s}^{\e}(\chi_F h^{i})\dif Y_{s}^{\e,i}\right|^{2}\Bigg]\no\\
&\leq& 5\e^{2}[\mE^{\mP^{\e}}(\Lambda_{T}^{\e})^{2}]^{\frac{1}{2}}\Bigg[\mE^{\mP^{\e}}\left|\int_{0}^{t}\rho_{s}^{\e}(\chi_F h^{i})\dif Y_{s}^{\e,i}\right|^{4}\Bigg]^{\frac{1}{2}}\no\\
&\leq&C\e^{2}\Bigg[\sum\limits_{i=1}^l\int_{0}^{T}\mE^{\mP^{\e}}|\rho_{s}^{\e}(\chi_F h^{i})|^{4}\dif s\Bigg]^{\frac{1}{2}}\leq C\e^{2}.
\label{j4}
\ee

Next, we deal with $J_{5}$. Set $\Psi(x,\mu,z):=F(x,\mu)h(x,\mu,z)-F(x,\mu)\bar{h}(x,\mu)$, and by the boundedness of $F$ and $(\mathbf{H}_{h})$, it holds that $\Psi(x,\mu,z)$ is bounded and Lipschitz continuous. Then, applying the H\"{o}lder inequality and the BDG inequality, one could obtain that
\be
J_{5}&=&5\mE\left|\int_{0}^{t}\rho_{s}^{\e}(\Psi^i)\dif Y_{s}^{\e,i}\right|^{2}
=5\mE^{\mP^{\e}}\left[\Lambda_{T}^{\e}\Big|\int_{0}^{t}\rho_{s}^{\e}(\Psi^i)\dif Y_{s}^{\e,i}\Big|^{2}\right]\no\\
&\leq& 5[\mE^{\mP^{\e}}(\Lambda_{T}^{\e})^{2}]^{\frac{1}{2}}\Bigg[\mE^{\mP^{\e}}\Big|\int_{0}^{t}\rho_{s}^{\e}(\Psi^i)\dif Y_{s}^{\e,i}\Big|^{4}\Bigg]^{\frac{1}{2}}
\leq C\Bigg[\sum\limits_{i=1}^l\mE^{\mP^{\e}}\int_{0}^{T}\Big|\rho_{s}^{\e}(\Psi^i)\Big|^{4}\dif s\Bigg]^{\frac{1}{2}}\no\\
&\leq&C\Bigg[\sum\limits_{i=1}^l\mE^{\mP^{\e}}\int_{0}^{T}\Big|\mE^{\mP^{\e}}[\Psi^i(X_{s}^{\e},\sL^{\mP}_{X_{s}^{\e}},Z_{s}^{\e})\Lambda_{s}^{\e}|\mathscr{F}_{s}^{Y^{\e}}]\Big|^{4}\dif s\Bigg]^{\frac{1}{2}}.
\label{j5es}
\ee

In the following, we are devoted to estimating the last term of the above inequality. First, we divide the time interval [0, T] into parts of size $\d$: $0=t_0<t_1<t_2<\cdots<t_N\leq T, N=[\frac{T}{\d}], t_{j+1}-t_{j}=\d, j=0,1,2,\cdots,N-1$, where $\d$ depends on $\e$. Then it holds that
\ce
&&\mE^{\mP^{\e}}\Bigg[\int_{0}^{T}\Big|\mE^{\mP^{\e}}[\Psi^i(X_{s}^{\e},\sL^{\mP}_{X_{s}^{\e}},Z_{s}^{\e})\Lambda_{s}^{\e}|\mathscr{F}_{s}^{Y^{\e}}]\Big|^{4}\dif s\Bigg]\\
&=&\sum\limits_{j=0}^{N-1}\mE^{\mP^{\e}}\Bigg[\int_{t_{j}}^{t_{j+1}}\Big|\mE^{\mP^{\e}}[\Psi^i(X_{s}^{\e},\sL^{\mP}_{X_{s}^{\e}},Z_{s}^{\e})\Lambda_{s}^{\e}|\mathscr{F}_{s}^{Y^{\e}}]\Big|^{4}\dif s\Bigg]\\
&&+\mE^{\mP^{\e}}\Bigg[\int_{N\d}^{T}\Big|\mE^{\mP^{\e}}[\Psi^i(X_{s}^{\e},\sL^{\mP}_{X_{s}^{\e}},Z_{s}^{\e})\Lambda_{s}^{\e}|\mathscr{F}_{s}^{Y^{\e}}]\Big|^{4}\dif s\Bigg].
\de
For clarity and simplicity, we use $[t,t+\d]$ instead of $[t_{j},t_{j+1}]$, and the analysis for the remainder term over the interval follows from the same argument. So, we know that
\be
&&\mE^{\mP^{\e}}\Bigg[\int_{0}^{T}\Big|\mE^{\mP^{\e}}[\Psi^i(X_{s}^{\e},\sL^{\mP}_{X_{s}^{\e}},Z_{s}^{\e})\Lambda_{s}^{\e}|\mathscr{F}_{s}^{Y^{\e}}]\Big|^{4}\dif s\Bigg]\no\\
&\leq&(N+1)\sup\limits_{t\in[0,T]}\Bigg[\mE^{\mP^{\e}}\int_{t}^{t+\d}\Big|\mE^{\mP^{\e}}[\Psi^i(X_{s}^{\e},\sL^{\mP}_{X_{s}^{\e}},Z_{s}^{\e})\Lambda_{s}^{\e}|\mathscr{F}_{s}^{Y^{\e}}]\Big|^{4}\dif s\Bigg]\no\\
&\leq&(\frac{T}{\d}+1)\sup\limits_{t\in[0,T]}\Bigg[\mE^{\mP^{\e}}\int_{t}^{t+\d}\Big|\mE^{\mP^{\e}}[\Psi^i(X_{s}^{\e},\sL^{\mP}_{X_{s}^{\e}},Z_{s}^{\e})\Lambda_{s}^{\e}|\mathscr{F}_{s}^{Y^{\e}}]\Big|^{4}\dif s\Bigg].
\label{ndt0}
\ee

Next, we study 
$$
\mE^{\mP^{\e}}\int_{t}^{t+\d}\Big|\mE^{\mP^{\e}}[\Psi^i(X_{s}^{\e},\sL^{\mP}_{X_{s}^{\e}},Z_{s}^{\e})\Lambda_{s}^{\e}|\mathscr{F}_{s}^{Y^{\e}}]\Big|^{4}\dif s.
$$
From the H\"{o}lder inequality and the Jensen inequality, it follows that
\be
&&\mE^{\mP^{\e}}\Bigg[\int_{t}^{t+\d}\Big|\mE^{\mP^{\e}}[\Psi^i(X_{s}^{\e},\sL^{\mP}_{X_{s}^{\e}},Z_{s}^{\e})\Lambda_{s}^{\e}|\mathscr{F}_{s}^{Y^{\e}}]\Big|^{4}\dif s\Bigg]\no\\
&=&\mE^{\mP^{\e}}\Bigg[\int_{t}^{t+\d}\Big|\mE^{\mP^{\e}}[\Psi^i(X_{s}^{\e},\sL^{\mP}_{X_{s}^{\e}},Z_{s}^{\e})\Lambda_{s}^{\e}
-\Psi^i(X_{t}^{\e},\sL^{\mP}_{X_{t}^{\e}},Z_{s}^{\e})\Lambda_{s}^{\e}\no\\
&&\quad+\Psi^i(X_{t}^{\e},\sL^{\mP}_{X_{t}^{\e}},Z_{s}^{\e})\Lambda_{s}^{\e}
-\Psi^i(X_{t}^{\e},\sL^{\mP}_{X_{t}^{\e}},Z_{s}^{\e})\Lambda_{t}^{\e}
+\Psi^i(X_{t}^{\e},\sL^{\mP}_{X_{t}^{\e}},Z_{s}^{\e})\Lambda_{t}^{\e}|\mathscr{F}_{s}^{Y^{\e}}]\Big|^{4}\dif s\Bigg]\no\\
&\leq&3^3\mE^{\mP^{\e}}\Bigg[\int_{t}^{t+\d}\Big|\mE^{\mP^{\e}}[\Psi^i(X_{s}^{\e},\sL^{\mP}_{X_{s}^{\e}},Z_{s}^{\e})\Lambda_{s}^{\e}
-\Psi^i(X_{t}^{\e},\sL^{\mP}_{X_{t}^{\e}},Z_{s}^{\e})\Lambda_{s}^{\e}|\mathscr{F}_{s}^{Y^{\e}}]\Big|^{4}\dif s\Bigg]\no\\
&&+3^3\mE^{\mP^{\e}}\Bigg[\int_{t}^{t+\d}\Big|\mE^{\mP^{\e}}[\Psi^i(X_{t}^{\e},\sL^{\mP}_{X_{t}^{\e}},Z_{s}^{\e})\Lambda_{s}^{\e}
-\Psi^i(X_{t}^{\e},\sL^{\mP}_{X_{t}^{\e}},Z_{s}^{\e})\Lambda_{t}^{\e}|\mathscr{F}_{s}^{Y^{\e}}]\Big|^{4}\dif s\Bigg]\no\\
&&+3^3\mE^{\mP^{\e}}\Bigg[\int_{t}^{t+\d}\Big|\mE^{\mP^{\e}}[\Psi^i(X_{t}^{\e},\sL^{\mP}_{X_{t}^{\e}},Z_{s}^{\e})\Lambda_{t}^{\e}|\mathscr{F}_{s}^{Y^{\e}}]\Big|^{4}\dif s\Bigg]\no\\
&\leq&3^3\int_{t}^{t+\d}\mE^{\mP^{\e}}\left|\left(\Psi^i(X_{s}^{\e},\sL^{\mP}_{X_{s}^{\e}},Z_{s}^{\e})
-\Psi^i(X_{t}^{\e},\sL^{\mP}_{X_{t}^{\e}},Z_{s}^{\e})\right)\Lambda_{s}^{\e}\right|^{4}\dif s\no\\
&&+3^3\int_{t}^{t+\d}\mE^{\mP^{\e}}\left|\Psi^i(X_{t}^{\e},\sL^{\mP}_{X_{t}^{\e}},Z_{s}^{\e})\left(\Lambda_{s}^{\e}-\Lambda_{t}^{\e}\right)\right|^{4}\dif s\no\\
&&+3^3\mE^{\mP^{\e}}\Bigg[\int_{t}^{t+\d}\Big|\mE^{\mP^{\e}}[\Psi^i(X_{t}^{\e},\sL^{\mP}_{X_{t}^{\e}},Z_{s}^{\e})\Lambda_{t}^{\e}|\mathscr{F}_{s}^{Y^{\e}}]\Big|^{4}\dif s\Bigg]\no\\
&=:&K_{1}+K_{2}+K_{3}.
\label{k1k2k3}
\ee

For $K_{1}$, by the H\"{o}lder inequality, and Lemma \ref{lambesti} and \ref{xediffp}, it holds that
\be
K_{1}&\leq&3^3\int_{t}^{t+\d}\(\mE^{\mP^{\e}}|\Psi^i(X_{s}^{\e},\sL^{\mP}_{X_{s}^{\e}},Z_{s}^{\e})-\Psi^i(X_{t}^{\e},\sL^{\mP}_{X_{t}^{\e}},Z_{s}^{\e})|^{8}\)^{\frac{1}{2}}[\mE^{\mP^{\e}}|
\Lambda_{s}^{\e}|^{8}]^{\frac{1}{2}}\dif s\no\\
&\leq& C\int_{t}^{t+\d}\(\mE^{\mP^{\e}}|\Psi^i(X_{s}^{\e},\sL^{\mP}_{X_{s}^{\e}},Z_{s}^{\e})-\Psi^i(X_{t}^{\e},\sL^{\mP}_{X_{t}^{\e}},Z_{s}^{\e})|^{8}\)^{\frac{1}{2}}\dif s\no\\
&\leq& C\int_{t}^{t+\d}\(\mE^{\mP^{\e}}\left(|X_{s}^{\e}-X_{t}^{\e}|^{8}+\mW_2(\sL^{\mP}_{X_{s}^{\e}},\sL^{\mP}_{X_{t}^{\e}})^{8}\right)\)^{\frac{1}{2}}\dif s\no\\
&\leq& C\int_{t}^{t+\d}\(\mE^{\mP^{\e}}|X_{s}^{\e}-X_{t}^{\e}|^{8}+\mE|X_{s}^{\e}-X_{t}^{\e}|^{8}\)^{\frac{1}{2}}\dif s\no\\
&=& C\int_{t}^{t+\d}\(\mE(\Lambda_{T}^{\e})^{-1}|X_{s}^{\e}-X_{t}^{\e}|^{8}+\mE|X_{s}^{\e}-X_{t}^{\e}|^{8}\)^{\frac{1}{2}}\dif s\no\\
&\leq& C\int_{t}^{t+\d}\((\mE|X_{s}^{\e}-X_{t}^{\e}|^{16})^{\frac{1}{2}}+\mE|X_{s}^{\e}-X_{t}^{\e}|^{8}\)^{\frac{1}{2}}\dif s\no\\
&\leq& C\int_{t}^{t+\d}(\d^{4}+\d^{2})\dif s\leq C(\d^{5}+\d^{3}),
\label{k1es}
\ee
where $\mW_2(\sL^{\mP}_{X_{s}^{\e}},\sL^{\mP}_{X_{t}^{\e}})^{2}\leq \mE|X_{s}^{\e}-X_{t}^{\e}|^{2}$.

From the boundedness of $\Psi$ and Lemma \ref{lambesti}, it follows that
\be
K_{2}\leq C\int_{t}^{t+\d}\mE^{\mP^{\e}}|\Lambda_{s}^{\e}-\Lambda_{t}^{\e}|^{4}\dif s\leq C\d^{3}.
\label{k2es}
\ee

To treat $K_{3}$, we introduce an auxiliary process as follows: for any $\mu\in\cP_2(\mR^m), z\in\mR^m$,
\be
\tilde{Z}_{s}^{\e,t,\mu,z}=z+\frac{1}{\e}\int_t^s b_{2}(\mu,\tilde{Z}_{r}^{\e,t,\mu,z})\dif r
+\frac{1}{\sqrt{\e}}\int_t^s\s_{2}(\mu,\tilde{Z}_{r}^{\e,t,\mu,z})\dif W_{r},\quad s\in[t, t+\d).
\label{hatz}
\ee
Then it holds that
\ce
K_{3}
&=&3^3\mE^{\mP^{\e}}\int_{t}^{t+\d}\Big|\mE^{\mP^{\e}}[\Psi^i(X_{t}^{\e},\sL^{\mP}_{X_{t}^{\e}},Z_{s}^{\e})\Lambda_{t}^{\e}-\Psi^i(X_{t}^{\e},\sL^{\mP}_{X_{t}^{\e}},\tilde{Z}_{s}^{\e,t,\sL^{\mP}_{X_{t}^{\e}},Z_{t}^{\e}})\Lambda_{t}^{\e}\\
&&\qquad +\Psi^i(X_{t}^{\e},\sL^{\mP}_{X_{t}^{\e}},\tilde{Z}_{s}^{\e,t,\sL^{\mP}_{X_{t}^{\e}},Z_{t}^{\e}})\Lambda_{t}^{\e}|\mathscr{F}_{s}^{Y^{\e}}]\Big|^{4}\dif s\\
&\leq&6^3\mE^{\mP^{\e}}\int_{t}^{t+\d}\Big|\mE^{\mP^{\e}}[\Psi^i(X_{t}^{\e},\sL^{\mP}_{X_{t}^{\e}},Z_{s}^{\e})\Lambda_{t}^{\e}-\Psi^i(X_{t}^{\e},\sL^{\mP}_{X_{t}^{\e}},\tilde{Z}_{s}^{\e,t,\sL^{\mP}_{X_{t}^{\e}},Z_{t}^{\e}})\Lambda_{t}^{\e}
|\mathscr{F}_{s}^{Y^{\e}}]\Big|^{4}\dif s\\
&&+6^3\mE^{\mP^{\e}}\int_{t}^{t+\d}\Big|\mE^{\mP^{\e}}[\Psi^i(X_{t}^{\e},\sL^{\mP}_{X_{t}^{\e}},\tilde{Z}_{s}^{\e,t,\sL^{\mP}_{X_{t}^{\e}},Z_{t}^{\e}})\Lambda_{t}^{\e}
|\mathscr{F}_{s}^{Y^{\e}}]\Big|^{4}\dif s\\
&=:&D_{1}+D_{2}.
\de

For $D_{1}$, based on the Jensen inequality, the H\"{o}lder inequality and (\ref{late}), it holds that
\be
D_{1}&\leq& C\int_{t}^{t+\d}\mE^{\mP^{\e}}|\Psi^i(X_{t}^{\e},\sL^{\mP}_{X_{t}^{\e}},Z_{s}^{\e})\Lambda_{t}^{\e}-\Psi^i(X_{t}^{\e},\sL^{\mP}_{X_{t}^{\e}},\tilde{Z}_{s}^{\e,t,\sL^{\mP}_{X_{t}^{\e}},Z_{t}^{\e}})\Lambda_{t}^{\e}
|^{4}\dif s\no\\
&\leq&C\int_{t}^{t+\d}\(\mE^{\mP^{\e}}|\Psi^i(X_{t}^{\e},\sL^{\mP}_{X_{t}^{\e}},Z_{s}^{\e})-\Psi^i(X_{t}^{\e},\sL^{\mP}_{X_{t}^{\e}},\tilde{Z}_{s}^{\e,t,\sL^{\mP}_{X_{t}^{\e}},Z_{t}^{\e}})|^{8}\)^{\frac{1}{2}}\(\mE^{\mP^{\e}}|\Lambda_{t}^{\e}
|^{8}\)^{\frac{1}{2}}\dif s\no\\
&\leq&C\int_{t}^{t+\d}\Big(\mE^{\mP^{\e}}|\Psi^i(X_{t}^{\e},\sL^{\mP}_{X_{t}^{\e}},Z_{s}^{\e})-\Psi^i(X_{t}^{\e},\sL^{\mP}_{X_{t}^{\e}},\tilde{Z}_{s}^{\e,t,\sL^{\mP}_{X_{t}^{\e}},Z_{t}^{\e}})|^{8}\Big)^{\frac{1}{2}}\dif s\no\\
&=&C\int_{t}^{t+\d}\Big(\mE^{\mP^{\e}}[\mE^{\mP^{\e}}[|\Psi^i(X_{t}^{\e},\sL^{\mP}_{X_{t}^{\e}},Z_{s}^{\e})-\Psi^i(X_{t}^{\e},\sL^{\mP}_{X_{t}^{\e}},\tilde{Z}_{s}^{\e,t,\sL^{\mP}_{X_{t}^{\e}},Z_{t}^{\e}})|^{8}|\sF_t^{X^{\e}}\vee\sF_t^{Z^{\e}}]]\Big)^{\frac{1}{2}}\dif s\no\\
&=&C\int_{t}^{t+\d}\Bigg(\mE^{\mP^{\e}}\Big[\mE^{\mP^{\e}}\Big[|\Psi^i(x,\sL^{\mP}_{X_{t}^{\e}},Z_{s}^{\e,t,z})-\Psi^i(x,\sL^{\mP}_{X_{t}^{\e}},\tilde{Z}_{s}^{\e,t,\sL^{\mP}_{X_{t}^{\e}},z})|^{8}\Big]\Big|_{(x,z)=(X_{t}^{\e},Z_{t}^{\e})}\Big]\Bigg)^{\frac{1}{2}}\dif s,\no\\
\label{did1}
\ee
where $\mathscr{F}_{t}^{X^{\e}}, \mathscr{F}_{t}^{Z^{\e}}$ denote the usual augmentation of $\sigma\{X_{r}^{\e},0\leq r \leq t\}, \sigma\{Z_{r}^{\e},0\leq r \leq t\}$, respectively, and$Z_{s}^{\e,t,z}$ solves the following equation: for $s\in[t,t+\d]$
\ce
Z_{s}^{\e,t,z}=z+\frac{1}{\e}\int_t^s b_{2}(\sL^{\mP}_{X_{r}^{\e}},Z_{r}^{\e,t,z})\dif r+\frac{1}{\sqrt{\e}}\int_t^s\s_{2}(\sL^{\mP}_{X_{r}^{\e}},Z_{r}^{\e,t,z})\dif W_{r}.
\de
Then Lemma \ref{-2es} implies that for $s\in[t,t+\d]$
\be
&&\mE^{\mP^{\e}}[|\Psi^i(x,\sL^{\mP}_{X_{t}^{\e}},Z_{s}^{\e,t,z})-\Psi^i(x,\sL^{\mP}_{X_{t}^{\e}},\tilde{Z}_{s}^{\e,t,\sL^{\mP}_{X_{t}^{\e}},z})|^{8}]\no\\
&\leq& C\mE^{\mP^{\e}}[|Z_{s}^{\e,t,z}-\tilde{Z}_{s}^{\e,t,\sL^{\mP}_{X_{t}^{\e}},z}|^{8}]\leq C\left(\mE[|Z_{s}^{\e,t,z}-\tilde{Z}_{s}^{\e,t,\sL^{\mP}_{X_{t}^{\e}},z}|^{16}]\right)^{\frac{1}{2}}(\mE[|\Lambda_{T}^{\e}|^{-2}])^{\frac{1}{2}}\no\\
&\leq& C\left(\mE[|Z_{s}^{\e,t,z}-\tilde{Z}_{s}^{\e,t,\sL^{\mP}_{X_{t}^{\e}},z}|^{16}]\right)^{\frac{1}{2}}\no\\
&\leq& C(1+|x_0|^{8}+|z_0|^{8})\Big((\frac{\d^{15}}{\e^{16}}+\frac{\d^{7}}{\e^{8}})(\d^{17}+\d^{9})e^{C(\frac{\d^{16}}{\e^{16}}+\frac{\d^{8}}{\e^{8}})}\Big)^{\frac{1}{2}},
\label{zexmuhatz}
\ee
where the last inequality is proved in the Appendix. Inserting the above inequality into (\ref{did1}), by Lemma \ref{xezebp} one can obtain that
\be
D_{1}\leq C\d\Big((\frac{\d^{15}}{\e^{16}}+\frac{\d^{7}}{\e^{8}})(\d^{17}+\d^{9})e^{C(\frac{\d^{16}}{\e^{16}}+\frac{\d^{8}}{\e^{8}})}\Big)^{\frac{1}{4}}.
\label{d1d11}
\ee

Next, we calculate $D_{2}$. Set
$$
\cH_{s}:=\mathscr{F}_{s}^{Y^{\e}}\vee\sF_t^{X^{\e}}\vee\sF_t^{Z^{\e}}, \quad s\in[t, t+\d),
$$
and from the tower property of the conditional expectation, the Jensen inequality, Lemma \ref{lambesti} and the H\"{o}lder inequality, it follows that
\ce
D_{2}&=&6^3\mE^{\mP^{\e}}\int_{t}^{t+\d}\Big|\mE^{\mP^{\e}}[\Psi^i(X_{t}^{\e},\sL^{\mP}_{X_{t}^{\e}},\tilde{Z}_{s}^{\e,t,\sL^{\mP}_{X_{t}^{\e}},Z_{t}^{\e}})\Lambda_{t}^{\e}
|\mathscr{F}_{s}^{Y^{\e}}]\Big|^{4}\dif s\no\\
&=&6^3\mE^{\mP^{\e}}\int_{t}^{t+\d}\Big|\mE^{\mP^{\e}}\[\mE^{\mP^{\e}}[\Psi^i(X_{t}^{\e},\sL^{\mP}_{X_{t}^{\e}},\tilde{Z}_{s}^{\e,t,\sL^{\mP}_{X_{t}^{\e}},Z_{t}^{\e}})\Lambda_{t}^{\e}
|\cH_{s}]\Big|\mathscr{F}_{s}^{Y^{\e}}\]\Big|^{4}\dif s\no\\
&=&6^3\mE^{\mP^{\e}}\int_{t}^{t+\d}\Big|\mE^{\mP^{\e}}\[\Lambda_{t}^{\e}\mE^{\mP^{\e}}[\Psi^i(X_{t}^{\e},\sL^{\mP}_{X_{t}^{\e}},\tilde{Z}_{s}^{\e,t,\sL^{\mP}_{X_{t}^{\e}},Z_{t}^{\e}})
|\cH_{s}]\Big|\mathscr{F}_{s}^{Y^{\e}}\]\Big|^{4}\dif s\no\\
&=&6^3\mE^{\mP^{\e}}\int_{t}^{t+\d}\Big|\mE^{\mP^{\e}}\[\Lambda_{t}^{\e}\mE^{\mP^{\e}}[\Psi^i(x,\sL^{\mP}_{X_{t}^{\e}},\tilde{Z}_{s}^{\e,t,\sL^{\mP}_{X_{t}^{\e}},z})]
|_{(x,z)=(X_{t}^{\e},Z_{t}^{\e})}\Big|\mathscr{F}_{s}^{Y^{\e}}\]\Big|^{4}\dif s\no\\
&\leq&6^3\int_{t}^{t+\d}\mE^{\mP^{\e}}\Big|\Lambda_{t}^{\e}\mE^{\mP^{\e}}[\Psi^i(x,\sL^{\mP}_{X_{t}^{\e}},\tilde{Z}_{s}^{\e,t,\sL^{\mP}_{X_{t}^{\e}},z})]
|_{(x,z)=(X_{t}^{\e},Z_{t}^{\e})}\Big|^{4}\dif s\no\\
&\leq&6^3\int_{t}^{t+\d}\(\mE^{\mP^{\e}}|\Lambda_{t}^{\e}|^{8}\)^{\frac{1}{2}}\Bigg(\mE^{\mP^{\e}}\Big|\mE^{\mP^{\e}}[\Psi^i(x,\sL^{\mP}_{X_{t}^{\e}},\tilde{Z}_{s}^{\e,t,\sL^{\mP}_{X_{t}^{\e}},z})]
|_{(x,z)=(X_{t}^{\e},Z_{t}^{\e})}\Big|^{8}\Bigg)^{\frac{1}{2}}\dif s\no\\
&\leq&C\int_{t}^{t+\d}\Bigg(\mE^{\mP^{\e}}\Big|\mE^{\mP^{\e}}[\Psi^i(x,\sL^{\mP}_{X_{t}^{\e}},\tilde{Z}_{s}^{\e,t,\sL^{\mP}_{X_{t}^{\e}},z})]
|_{(x,z)=(X_{t}^{\e},Z_{t}^{\e})}\Big|^{8}\Bigg)^{\frac{1}{2}}\dif s\no\\
&\leq&C\d^{\frac{1}{2}}\Bigg(\int_{t}^{t+\d}\mE^{\mP^{\e}}\Big|\mE^{\mP^{\e}}[\Psi^i(x,\sL^{\mP}_{X_{t}^{\e}},\tilde{Z}_{s}^{\e,t,\sL^{\mP}_{X_{t}^{\e}},z})]
|_{(x,z)=(X_{t}^{\e},Z_{t}^{\e})}\Big|^{8}\dif s\Bigg)^{\frac{1}{2}}.
\de

Next, on one hand, it holds that
\ce
\tilde{Z}^{\e,t,\mu, z}_{\e s+t}&=&z+\frac{1}{\e} \int_{t}^{\e s+t}b_2(\mu, \tilde{Z}^{\e,t,\mu, z}_{r})\dif r+\frac{1}{\sqrt{\e}} \int_{t}^{\e s+t} \s_2(\mu, \tilde{Z}^{\e, t, \mu, z}_{r})\dif W_r\\
&=&z+\frac{1}{\e} \int_{0}^{\e s}b_2(\mu, \tilde{Z}^{\e,\mu, z}_{u+t})\dif u+\frac{1}{\sqrt{\e}} \int_{0}^{\e s} \s_2(\mu, \tilde{Z}^{\e,t,\mu, z}_{u+t})\dif \tilde{W}_u\\
&=&z+\int_{0}^{s}b_2(\mu, \tilde{Z}^{\e, t, \mu, z}_{\e v+t})\dif v+\int_{0}^{s} \s_2(\mu, \tilde{Z}^{\e, t, \mu, z}_{\e v+t})\dif \check{\tilde{W}}_v,
\de
where $\tilde{W}_u:=W_{u+k\d}-W_{k\d}$ and $\check{\tilde{W}}_v:=\frac{1}{\sqrt{\e}}\tilde{W}_{\e v}$ are two $m$-dimensional standard Brownian motions. On the other hand, note that the frozen equation (\ref{Eq2}) is written as
\ce
Z_{s}^{\mu, z}=z+\int_0^s b_{2}(\mu,Z_{r}^{\mu, z})\dif r+\int_0^s\s_{2}(\mu,Z_{r}^{\mu, z})\dif W_{r}.
\de
Thus, for $s\in[0,\d/\e]$, $\tilde{Z}^{\e, t, \mu, z}_{\e s+t}$ and $Z_{s}^{\mu, z}$ have the same distribution, which implies that
\be
D_{2}&\leq&C\d^{\frac{1}{2}}\Bigg(\int_{t}^{t+\d}\mE^{\mP^{\e}}\Big[P_{\frac{s-t}{\e}}^{\sL^{\mP}_{X_{t}^{\e}}}\Psi^i(X_{t}^{\e},\sL^{\mP}_{X_{t}^{\e}},Z_{t}^{\e})
\Big]^{8}\dif s\Bigg)^{\frac{1}{2}}\no\\
&=&C\d^{\frac{1}{2}}\e^{\frac{1}{2}}\Bigg(\int_{0}^{\frac{\d}{\e}}\mE^{\mP^{\e}}\Big[P_{r}^{\sL^{\mP}_{X_{t}^{\e}}}\Psi^i(X_{t}^{\e},\sL^{\mP}_{X_{t}^{\e}},Z_{t}^{\e})
\Big]^{8}\dif r\Bigg)^{\frac{1}{2}}\no\\
&\leq&C\d^{\frac{1}{2}}\e^{\frac{1}{2}}\Bigg(\int_{0}^{\infty}\mE^{\mP^{\e}}\Big[P_{r}^{\sL^{\mP}_{X_{t}^{\e}}}\Psi^i(X_{t}^{\e},\sL^{\mP}_{X_{t}^{\e}},Z_{t}^{\e})
\Big]^{8}\dif r\Bigg)^{\frac{1}{2}}.
\label{d2d2}
\ee

In the following, we deal with  $|P_{r}^{\sL^{\mP}_{X_{t}^{\e}}}\Psi^i(X_{t}^{\e},\sL^{\mP}_{X_{t}^{\e}},Z_{t}^{\e})|^8$. Set $\bar{\Psi}(x,\mu):=\int_{\mR^{m}}\Psi(x,\mu,z)\nu^{\mu}(\dif z)$, and by the same deduction as that of (\ref{semiesti}) it holds that 
\be
&&\Big|P_{r}^{\sL^{\mP}_{X_{t}^{\e}}}\Psi^i(X_{t}^{\e},\sL^{\mP}_{X_{t}^{\e}},Z_{t}^{\e})-\bar{\Psi}^i(X_{t}^{\e},\sL^{\mP}_{X_{t}^{\e}})\Big|^{2}\no\\
&=&\Big|\mE\left[\Psi^i(x,\mu,Z_{r}^{\mu,z})-\bar{\Psi}^i(x,\mu)\right]\Big|_{(x,\mu,z)=(X_{t}^{\e},\sL^{\mP}_{X_{t}^{\e}},Z_{t}^{\e})}\Big|^{2}\no\\
&\leq& Ce^{-\a'r}(1+\|\sL^{\mP}_{X_{t}^{\e}}\|^{2}+|Z_{t}^{\e}|^{2}).
\label{ppsi}
\ee
Note that $\bar{\Psi}(x,\mu)=0$. Thus, one can obtain that
\ce
\Big|P_{r}^{\sL^{\mP}_{X_{t}^{\e}}}\Psi^i(X_{t}^{\e},\sL^{\mP}_{X_{t}^{\e}},Z_{t}^{\e})\Big|^{2}
\leq Ce^{-\a' r}(1+\|\sL^{\mP}_{X_{t}^{\e}}\|^{2}+|Z_{t}^{\e}|^{2}),
\de
and 
\ce
\Big|P_{r}^{\sL^{\mP}_{X_{t}^{\e}}}\Psi^i(X_{t}^{\e},\sL^{\mP}_{X_{t}^{\e}},Z_{t}^{\e})\Big|^{8}
\leq Ce^{-4\a' r}(1+\|\sL^{\mP}_{X_{t}^{\e}}\|^{8}+|Z_{t}^{\e}|^{8}).
\de

Next, by inserting the above inequality in (\ref{d2d2}), it holds that
\be
D_{2}&\leq& C\d^{\frac{1}{2}}\e^{\frac{1}{2}}\Bigg(\int_{0}^{\infty}e^{-4\a' r}\mE^{\mP^{\e}}(1+\|\sL^{\mP}_{X_{t}^{\e}}\|^{8}+|Z_{t}^{\e}|^{8})\dif r\Bigg)^{\frac{1}{2}}\no\\
&\leq&C\d^{\frac{1}{2}}\e^{\frac{1}{2}}\left(\mE^{\mP^{\e}}(1+\|\sL^{\mP}_{X_{t}^{\e}}\|^{8}+|Z_{t}^{\e}|^{8})\right)^{\frac{1}{2}}\no\\
&\leq&C\d^{\frac{1}{2}}\e^{\frac{1}{2}}\left(\mE(1+\|\sL^{\mP}_{X_{t}^{\e}}\|^{16}+|Z_{t}^{\e}|^{16})\right)^{\frac{1}{4}}\no\\
&\leq&C\d^{\frac{1}{2}}\e^{\frac{1}{2}},
\label{d2boun}
\ee
where Lemma \ref{xezebp} is used in the last inequality. 

So, combining (\ref{d2boun}) with (\ref{d1d11}), we know that
\be
K_{3}\leq C\d\Big((\frac{\d^{15}}{\e^{16}}+\frac{\d^{7}}{\e^{8}})(\d^{17}+\d^{9})e^{C(\frac{\d^{16}}{\e^{16}}+\frac{\d^{8}}{\e^{8}})}\Big)^{\frac{1}{4}}
+C\d^{\frac{1}{2}}\e^{\frac{1}{2}}.
\label{k3es}
\ee
Inserting (\ref{k1es}), (\ref{k2es}) and (\ref{k3es}) into (\ref{k1k2k3}), one can get that
\ce
&&\mE^{\mP^{\e}}\Bigg[\int_{t}^{t+\d}\Big|\mE^{\mP^{\e}}[\Psi^i(X_{s}^{\e},\sL^{\mP}_{X_{s}^{\e}},Z_{s}^{\e})\Lambda_{s}^{\e}|\mathscr{F}_{s}^{Y^{\e}}]\Big|^{4}\dif s\Bigg]\\
&\leq& C(\d^{5}+\d^{3})+C\d^{3}+ C\d\Big((\frac{\d^{15}}{\e^{16}}+\frac{\d^{7}}{\e^{8}})(\d^{17}+\d^{9})e^{C(\frac{\d^{16}}{\e^{16}}+\frac{\d^{8}}{\e^{8}})}\Big)^{\frac{1}{4}}
+C\d^{\frac{1}{2}}\e^{\frac{1}{2}},
\de
which together with (\ref{ndt0}) implies that
\ce
&&\mE^{\mP^{\e}}\Bigg[\int_{0}^{T}\Big|\mE^{\mP^{\e}}[\Psi^i(X_{s}^{\e},\sL^{\mP}_{X_{s}^{\e}},Z_{s}^{\e})\Lambda_{s}^{\e}|\mathscr{F}_{s}^{Y^{\e}}]\Big|^{4}\dif s\Bigg]\\
&\leq&(\frac{T}{\d}+1) \[C(\d^{5}+\d^{3})+C\d^{3}+ C\d\Big((\frac{\d^{15}}{\e^{16}}+\frac{\d^{7}}{\e^{8}})(\d^{17}+\d^{9})e^{C(\frac{\d^{16}}{\e^{16}}+\frac{\d^{8}}{\e^{8}})}\Big)^{\frac{1}{4}}
+C\d^{\frac{1}{2}}\e^{\frac{1}{2}}\]\\
&\leq&C(T+\d) \[(\d^{4}+\d^{2})+\d^{2}+\Big((\frac{\d^{15}}{\e^{16}}+\frac{\d^{7}}{\e^{8}})(\d^{17}+\d^{9})e^{C(\frac{\d^{16}}{\e^{16}}+\frac{\d^{8}}{\e^{8}})}\Big)^{\frac{1}{4}}
+\frac{\e^{\frac{1}{2}}}{\d^{\frac{1}{2}}}\].
\de
From the above inequality and (\ref{j5es}), it follows that
\be
J_{5}\leq C(T+\d)^{\frac{1}{2}} \[(\d^{4}+\d^{2})+\d^{2}+\Big((\frac{\d^{15}}{\e^{16}}+\frac{\d^{7}}{\e^{8}})(\d^{17}+\d^{9})e^{C(\frac{\d^{16}}{\e^{16}}+\frac{\d^{8}}{\e^{8}})}\Big)^{\frac{1}{4}}
+\frac{\e^{\frac{1}{2}}}{\d^{\frac{1}{2}}}\]^{\frac{1}{2}}.
\label{j5}
\ee

Finally, inserting (\ref{j1})-(\ref{j4}) and (\ref{j5}) into (\ref{mxifc}), we have
\ce
&&\mE\left|\xi_{t}^{\e}(F)-\int_{0}^{t}\xi_{s}^{\e}(\bar{\cL}F)\dif s-\int_{0}^{t}\xi_{s}^{\e}(F \bar{h}^{i})\dif Y_{s}^{\e,i}\right|^{2}\\
&\leq& C\e^{2}+C(T+\d)^{\frac{1}{2}} \[(\d^{4}+\d^{2})+\d^{2}+\Big((\frac{\d^{15}}{\e^{16}}+\frac{\d^{7}}{\e^{8}})(\d^{17}+\d^{9})e^{C(\frac{\d^{16}}{\e^{16}}+\frac{\d^{8}}{\e^{8}})}\Big)^{\frac{1}{4}}
+\frac{\e^{\frac{1}{2}}}{\d^{\frac{1}{2}}}\]^{\frac{1}{2}}.
\de
By choosing $\d=\e(-ln\e)^{\frac{1}{32}}$ and taking the limit on both sides of the above inequality, it holds that
\ce
\lim_{\e\rightarrow0}\mE\left|\xi_{t}^{\e}(F)-\int_{0}^{t}\xi_{s}^{\e}(\bar{\cL}F)\dif s-\int_{0}^{t}\xi_{s}^{\e}(F \bar{h}^{i})\dif Y_{s}^{\e,i}\right|^{2}=0.
\de

{\bf Step 2.} We prove that there exists a subsequence $\{\xi^{\e_k}\}$ which converges weakly to $0$.

By Lemma \ref{corbr} and {\bf Step 1}, we know that there exists a subsequence $\{\xi^{\e_k}\}$ which converges weakly to $\xi$ in $C([0,T], \cM(\mR^{n}\times\sP_{2}(\mR^{n}))$ as $k\rightarrow\infty$, and furthermore for $t\in[0,T]$ and $F\in \mC^{4,(2,2)}_{b}(\mR^{n}\times\sP_{2}(\mR^{n}))$, $\xi_{t}(F)$ satisfies the equation
\be
\xi_{t}(F)-\int_{0}^{t}\xi_{s}(\bar{\cL}F)\dif s-\int_{0}^{t}\xi_{s}(F \bar{h}^{i})\dif Y_{s}=0,
\label{xiequa}
\ee
where $Y$ is a $l$-dimensional Brownian motion. For Eq.(\ref{xiequa}), by \cite[Theorem 4.9]{MQ}, it holds that its solutions are unique. Besides, note that $0$ is a solution to Eq.(\ref{xiequa}). Thus, $\xi=0$, that is, $\{\xi^{\e_k}\}$ converges weakly to $0$. The proof is complete.
\end{proof}

Now, it is the position to prove Theorem \ref{convfilt}.

{\bf Proof of Theorem \ref{convfilt}.}
For $t\in[0,T]$ and $F\in \mC^{4,(2,2)}_{b}(\mR^{n}\times\sP_{2}(\mR^{n}))$, it holds that 
\ce
\pi_{t}^{\e,x,\mu}(F)-\bar{\pi}_{t}(F)
=\frac{\rho_{t}^{\e,x,\mu}(F)-\bar{\rho}_{t}(F)}{\bar{\rho}_{t}(1)}
-\pi_{t}^{\e,x,\mu}(F)\frac{\rho_{t}^{\e,x,\mu}(1)-\bar{\rho}_{t}(1)}{\bar{\rho}_{t}(1)}=\frac{\xi_t^\e(F)}{\bar{\rho}_{t}(1)}
-\pi_{t}^{\e,x,\mu}(F)\frac{\xi_t^\e(1)}{\bar{\rho}_{t}(1)}.
\de
So, based on Lemma \ref{rhbi} and \ref{corbr2}, and the boundedness of $\pi_{t}^{\e,x,\mu}(F)$, we obtain that $\pi_{t}^{\e,x,\mu}$ converges weakly to $\bar{\pi}_{t}$ as $\e\rightarrow0$. The proof is complete.

\section{An example}\label{exam}

Now let us present an example to explain our results. 

\bx
Consider the following multiscale McKean-Vlasov stochastic system on $\mR^{n} \times \mR^{m}$:
\be\left\{\begin{array}{l}
\dif X_{t}^{\e}=b_{1}(X_{t}^{\e}, \sL^{\mP}_{X_{t}^{\e}}, Z_{t}^{\e})\dif t+\dif B_{t},\\
X_{0}^{\e}=x_0,\quad  0\leq t\leq T,\\
\dif Z_{t}^{\e}=\frac{1}{\e}b_{2}(\sL^{\mP}_{X_{t}^{\e}},Z_{t}^{\e})\dif t+\frac{1}{\sqrt{\e}}\dif W_{t},\\
Z_{0}^{\e}=z_0,\quad  0\leq t\leq T,
\end{array}
\right.
\label{slfaex}
\ee
where $b_{1}(x,\mu,z)=\int_{\mR^n}\tilde{b}_{1}(x+u,z)\mu(\dif u)$, $b_{2}(\mu,z)=\int_{\mR^n}\tilde{b}_{2}(u,z)\mu(\dif u)$ and $\tilde{b}_{1}: \mR^n\times\mR^m\mapsto\mR^n$, $\tilde{b}_{2}: \mR^n\times\mR^m\mapsto\mR^m$ are Borel measurable.

Assume:

$(i)$ $\partial_x \tilde{b}_{1}(x, z), \partial_z \tilde{b}_{1}(x, z), \partial_{xx} \tilde{b}_{1}(x, z), \partial_{xz} \tilde{b}_{1}(x, z)$ exist for any $x\in\mR^n, z\in\mR^m$. Moreover, all these partial derivatives are uniformly bounded and Lipschitz continuous w.r.t. $z$ uniformly in $x$.

$(ii)$ There exists $\beta>0$ such that for any $x \in \mathbb{R}^n$ and $z_1, z_2 \in \mathbb{R}^m$,
$$
2\left\langle  z_1-z_2, \tilde{b}_{2}\left(x, z_1\right)-\tilde{b}_{2}\left(x, z_2\right)\right\rangle \leqslant-\beta\left|z_1-z_2\right|^2
$$

$(iii)$ $\partial_x \tilde{b}_{2}(x, z), \partial_z \tilde{b}_{2}(x, z), \partial_{xx} \tilde{b}_{2}(x, z), \partial_{xz} \tilde{b}_{2}(x, z)$ exist for any $x\in\mR^n, z\in\mR^m$. Moreover, all these partial derivatives are uniformly bounded and Lipschitz continuous w.r.t. $z$ uniformly in $x$.

Note that
\ce
&&\partial_\mu b_{1}(x,\mu,z)(y)=\partial_x \tilde{b}_{1}(x+y, z), \quad \partial_y\partial_\mu b_{1}(x,\mu,z)(y)=\partial_{xx} \tilde{b}_{1}(x+y, z),\\
&&\partial_\mu b_{2}(\mu,z)(y)=\partial_x \tilde{b}_{2}(y, z), \quad \partial_y\partial_\mu b_{2}(\mu,z)(y)=\partial_{xx} \tilde{b}_{2}(y, z).
\de
So, if $(i)$ $(ii)$ $(iii)$ and $\b>5L_{b_{2}, \s_{2}}$ hold, where $L_{b_{2}, \s_{2}}:=2(\max\{\|\partial_x \tilde{b}_{2}\|, \|\partial_z \tilde{b}_{2}\|\})^2$, $b_1,\s_1,b_2,\s_2$ satisfy $(\mathbf{H}^1_{b_{1}, \s_{1}})$ $(\mathbf{H}^2_{\s_{1}})$ $(\mathbf{H}^3_{b_{1}, \s_{1}})$ $(\mathbf{H}^1_{b_{2}, \s_{2}})$-$(\mathbf{H}^3_{b_{2}, \s_{2}})$.Therefore, by Theorem \ref{weakconv}, we know that $\{X_{t}^{\e}, t\in[0,T]\}$ converges weakly to $\{\bar{X}_{t}, t\in[0,T]\}$ in $C([0,T],\mR^{n})$, where $\bar{X}$ solves the corresponding averaged equation.

Next, given an observation process $Y_{t}^{\e}$, i.e.
\be
Y_{t}^{\e}=V_{t}+\int_{0}^{t}h(X_{s}^{\e},\sL^{\mP}_{X_{s}^{\e}},Z_{s}^{\e})\dif s,
\label{byteex}
\ee
where $h(x,\mu,z)=\int_{\mR^n}\sin|x+u|\mu(\dif u)+\sin|z|$. Then it is easy to justify that $h$ satisfies $(\mathbf{H}_{h})$. If $(i)$ $(ii)$ $(iii)$ and $\b>(2p+1)L_{b_{2}, \s_{2}}$ hold, $b_1,\s_1,b_2,\s_2,h$ satisfy $(\mathbf{H}^1_{b_{1}, \s_{1}})$ $(\mathbf{H}^2_{\s_{1}})$ $(\mathbf{H}^3_{b_{1}, \s_{1}})$ $(\mathbf{H}^{1}_{b_{2}, \s_{2}})$ $(\mathbf{H}^{2'}_{b_{2}, \s_{2}})$ $(\mathbf{H}^{3}_{b_{2}, \s_{2}})$ and $(\mathbf{H}_{h})$. Hence, by Theorem \ref{convfilt}, we have that for any $t\in[0,T]$, $\pi_{t}^{\e,x,\mu}$ converges weakly to $\bar{\pi}_{t}$ as $\e\rightarrow 0$.
\ex

\section{Appendix}\label{app}

In this section, we prove (\ref{zexmuhatz}). 

{\bf Proof of (\ref{zexmuhatz}).} By the H\"{o}lder inequality, the BDG inequality, $(\mathbf{H}^1_{b_{2}, \s_{2}})$ and Lemma \ref{xediffp}, we get that for $s\in[t,t+\d]$
\ce
&&\mE[|Z_{s}^{\e,t,z}-\tilde{Z}_{s}^{\e,t,\sL^{\mP}_{X_{t}^{\e}},z}|^{16}]\\
&=&\mE\Bigg[\Big|\frac{1}{\e}\int_{t}^{s}\(b_{2}(\sL^{\mP}_{X^\e_{r}},Z_{r}^{\e,t,z})-b_{2}(\sL^{\mP}_{X_{t}^{\e}},\tilde{Z}_{r}^{\e,t,\sL^{\mP}_{X_{t}^{\e}},z})\)\dif r\\
&&\quad+\frac{1}{\sqrt{\e}}\int_{t}^{s}\(\s_{2}(\sL^{\mP}_{X^\e_{r}},Z_{r}^{\e,t,z})-\s_{2}(\sL^{\mP}_{X_{t}^{\e}},\tilde{Z}_{r}^{\e,t,\sL^{\mP}_{X_{t}^{\e}},z})\)\dif W_{r}\Big|^{16}\Bigg]\\
&\leq& \frac{2^{15}}{\e^{16}}\mE\Bigg[\Big|\int_{t}^{s}\(b_{2}(\sL^{\mP}_{X^\e_{r}},Z_{r}^{\e,t,z})-b_{2}(\sL^{\mP}_{X_{t}^{\e}},\tilde{Z}_{r}^{\e,t,\sL^{\mP}_{X_{t}^{\e}},z})\)\dif r\Big|^{16}\Bigg]\\
&&+\frac{2^{15}}{\e^{8}}\mE\Bigg[\Big|\int_{t}^{s}\(\s_{2}(\sL^{\mP}_{X^\e_{r}},Z_{r}^{\e,t,z})-\s_{2}(\sL^{\mP}_{X_{t}^{\e}},\tilde{Z}_{r}^{\e,t,\sL^{\mP}_{X_{t}^{\e}},z})\)\dif W_{r}\Big|^{16}\Bigg]\\
&\leq& C\frac{\d^{15}}{\e^{16}}\int_{t}^{s}\mE\Big|b_{2}(\sL^{\mP}_{X^\e_{r}},Z_{r}^{\e,t,z})-b_{2}(\sL^{\mP}_{X_{t}^{\e}},\tilde{Z}_{r}^{\e,t,\sL^{\mP}_{X_{t}^{\e}},z})\Big|^{16}\dif r\\
&&+C\frac{\d^{7}}{\e^{8}}\int_{t}^{s}\mE\|\s_{2}(\sL^{\mP}_{X^\e_{r}},Z_{r}^{\e,t,z})-\s_{2}(\sL^{\mP}_{X_{t}^{\e}},\tilde{Z}_{r}^{\e,t,\sL^{\mP}_{X_{t}^{\e}},z})\|^{16}\dif r\\
&=& C\frac{\d^{15}}{\e^{16}}\int_{t}^{s}\mE\Big|b_{2}(\sL^{\mP}_{X^\e_{r}},Z_{r}^{\e,t,z})-b_{2}(\sL^{\mP}_{X^\e_{t}},Z_{r}^{\e,t,z})+b_{2}(\sL^{\mP}_{X_{t}^{\e}},Z_{r}^{\e,t,z})-b_{2}(\sL^{\mP}_{X_{t}^{\e}},\tilde{Z}_{r}^{\e,t,\sL^{\mP}_{X_{t}^{\e}},z})\Big|^{16}\dif r\\
&&+C\frac{\d^{7}}{\e^{8}}\int_{t}^{s}\mE\|\s_{2}(\sL^{\mP}_{X^\e_{r}},Z_{r}^{\e,t,z})-\s_{2}(\sL^{\mP}_{X_{t}^{\e}},Z_{r}^{\e,t,z})+\s_{2}(\sL^{\mP}_{X_{t}^{\e}},Z_{r}^{\e,t,z})-\s_{2}(\sL^{\mP}_{X_{t}^{\e}},\tilde{Z}_{r}^{\e,t,\sL^{\mP}_{X_{t}^{\e}},z})\|^{16}\dif r\\
&\leq& C\frac{\d^{15}}{\e^{16}}\int_{t}^{s}\mE\Big|b_{2}(\sL^{\mP}_{X^\e_{r}},Z_{r}^{\e,t,z})-b_{2}(\sL^{\mP}_{X_{t}^{\e}},Z_{r}^{\e,t,z})\Big|^{16}\dif r\\
&&+C\frac{\d^{15}}{\e^{16}}\int_{t}^{s}\mE\Big|b_{2}(\sL^{\mP}_{X_{t}^{\e}},Z_{r}^{\e,t,z})-b_{2}(\sL^{\mP}_{X_{t}^{\e}},\tilde{Z}_{r}^{\e,t,\sL^{\mP}_{X_{t}^{\e}},z})\Big|^{16}\dif r\\
&&+C\frac{\d^{7}}{\e^{8}}\int_{t}^{s}\mE\|\s_{2}(\sL^{\mP}_{X^\e_{r}},Z_{r}^{\e,t,z})-\s_{2}(\sL^{\mP}_{X_{t}^{\e}},Z_{r}^{\e,\sL^{\mP}_{X_{t}^{\e}},z})\|^{16}\dif r\\
&&+C\frac{\d^{7}}{\e^{8}}\int_{t}^{s}\mE\|\s_{2}(\sL^{\mP}_{X_{t}^{\e}},Z_{r}^{\e,t,z})-\s_{2}(\sL^{\mP}_{X_{t}^{\e}},\tilde{Z}_{r}^{\e,t,\sL^{\mP}_{X_{t}^{\e}},z})\|^{16}\dif r\\
&\leq& C(\frac{\d^{15}}{\e^{16}}+\frac{\d^{7}}{\e^{8}})\int_{t}^{s}\mE L^8_{b_{2},\s_{2}}\mW_2^{16}(\sL^{\mP}_{X_{r}^{\e}},\sL^{\mP}_{X_{t}^{\e}})\dif r+C(\frac{\d^{15}}{\e^{16}}+\frac{\d^{7}}{\e^{8}})\int_{t}^{s}\mE[L_{b_{2},\s_{2}}|Z_{r}^{\e,t,z}-\tilde{Z}_{r}^{\e,t,\sL^{\mP}_{X_{t}^{\e}},z}|^{2}]^{8}\dif r\\
&\leq& C(\frac{\d^{15}}{\e^{16}}+\frac{\d^{7}}{\e^{8}})\int_{t}^{s}\mE|X_{r}^{\e}-X_{t}^{\e}|^{16}\dif r+C(\frac{\d^{15}}{\e^{16}}+\frac{\d^{7}}{\e^{8}})\int_{t}^{s}\mE|Z_{r}^{\e,t,z}-\tilde{Z}_{r}^{\e,t,\sL^{\mP}_{X_{t}^{\e}},z}|^{16}\dif r\\
&\leq& C(1+|x_0|^{16}+|z_0|^{16})(\frac{\d^{15}}{\e^{16}}+\frac{\d^{7}}{\e^{8}})\int_{t}^{s}(\d^{16}+\d^{8})\dif r+C(\frac{\d^{15}}{\e^{16}}+\frac{\d^{7}}{\e^{8}})\int_{t}^{s}\mE|Z_{r}^{\e,t,z}-\tilde{Z}_{r}^{\e,t,\sL^{\mP}_{X_{t}^{\e}},z}|^{16}\dif r\\
&\leq&C(1+|x_0|^{16}+|z_0|^{16})(\frac{\d^{15}}{\e^{16}}+\frac{\d^{7}}{\e^{8}})(\d^{17}+\d^{9})
+C(\frac{\d^{15}}{\e^{16}}+\frac{\d^{7}}{\e^{8}})\int_{t}^{s}\mE|Z_{r}^{\e,t,z}-\tilde{Z}_{r}^{\e,t,\sL^{\mP}_{X_{t}^{\e}},z}|^{16}\dif r.
\de
Then the Gronwall inequality implies that
\ce
\mE[|Z_{s}^{\e,t,z}-\tilde{Z}_{s}^{\e,t,\sL^{\mP}_{X_{t}^{\e}},z}|^{16}]\leq C(1+|x_0|^{16}+|z_0|^{16})(\frac{\d^{15}}{\e^{16}}+\frac{\d^{7}}{\e^{8}})(\d^{17}+\d^{9})e^{C(\frac{\d^{16}}{\e^{16}}+\frac{\d^{8}}{\e^{8}})}.
\de

\end{document}